\documentclass[10pt, reqno]{article}

\usepackage{helvet}
\usepackage{amsfonts}
\usepackage{latexsym}
\usepackage{amsmath}
\usepackage{amsthm,amssymb,mathrsfs,amstext}
\usepackage{graphicx}
\input epsf.tex
\usepackage{epsf}
\usepackage{epsfig}
\usepackage[all]{xy}
\font\msbm=msbm10

\numberwithin{equation}{section}
 \def\cE{\mathcal{E}}
\def\cF{\mathcal{F}}              

 \def\cH{\mathcal{H}}
 \def\cB{\mathcal{B}}
 \def\cG{\mathcal{G}}
 \def\cM{\mathcal{M}}
 \def\cK{\mathcal{K}}
 \def\cU{\mathcal{U}}
 \def\cA{\mathcal{A}}

 \def\cT{\mathcal{T}}
\def\cV{\mathcal{V}}

\theoremstyle{plain}
\newtheorem{Theorem}{Theorem}[section]
\newtheorem{lemma}[Theorem]{Lemma}
\newtheorem{corollary}[Theorem]{Corollary}
\newtheorem{proposition}[Theorem]{Proposition}

\theoremstyle{definition}
\newtheorem{example}{Example}[section]
\newtheorem{examples}[example]{Example}
\newtheorem{definition}{Definition}[section]
\theoremstyle{remark}
\newtheorem{remark}{Remark}[section]
\title{Continuous Frames, Function Spaces, and the Discretization Problem}
\author{Massimo Fornasier\footnote{The author acknowledges the partial support of the Intra-European Individual Marie Curie Fellowship, project FTFDORF-FP6-501018, and the hospitality of NuHAG (Numerical Harmonic Analysis Group), Fakult\"at f\"ur Mathematik, Universit\"at Wien, Austria, during the preparation of this work.}, Holger Rauhut}
\def\mathbb#1{\hbox{\msbm{#1}}}
\newcommand{\N}{{\mathbb{N}}}
\newcommand{\R}{{\mathbb{R}}}

\newcommand{\C}{{\mathbb{C}}}

\newcommand{\A}{{\cal{A}}}

\newcommand{\G}{{\cal{G}}}
\newcommand{\K}{{\cal{K}}}
\renewcommand{\H}{{\cal{H}}}
\newcommand{\Co}{{\mathsf{Co}}}
\newcommand{\TCo}{\widetilde{{\mathsf{Co}}}}
\newcommand{\for}{{\mbox{ for all }}}
\newcommand{\spann}{\operatorname{span}}
\newcommand{\osc}{\operatorname{osc}}
\newcommand{\on}{{|\!|\!|}}
\newcommand{\UP}{U_{\Phi}}
\newcommand{\SP}{S_{\Phi}}

\newcommand{\beq}{\begin{eqnarray}}
\newcommand{\eeq}{\end{eqnarray}}
\newcommand{\beqn}{\begin{eqnarray*}}
\newcommand{\eeqn}{\end{eqnarray*}}
\newcommand{\id}{\operatorname{Id}}
\newcommand{\F}{{\cal F}}
\newcommand{\ol}{\overline}
\newcommand{\supp}{\operatorname{supp}}
\newcommand{\esssup}{\operatorname{ess\,sup}}
\newcommand{\bdiscrete}{\flat}
\newcommand{\ddiscrete}{\natural}
\newcommand{\admissible}{admissible }
\newcommand{\moderate}{moderate admissible }

\renewcommand{\qed}{\rule{2.5mm}{2.5mm}}
\newenvironment{Proof}{\noindent
{\bf\underline{Proof:} }}
{\hspace*{\fill}\qed\vskip1em}

\begin{document}
\maketitle
\begin{abstract}
A continuous frame is a family of vectors in a 
Hilbert space which allows 
reproductions of arbitrary elements by continuous superpositions.
Associated to a given continuous frame we construct certain Banach spaces.
Many classical function spaces can be identified as such spaces.
We provide a general method to derive Banach frames and atomic decompositions
for these Banach spaces by sampling the continuous frame. This is done by generalizing
the coorbit space theory developed by Feichtinger and Gr\"ochenig.
As an important tool the concept of 
localization of frames 
is extended to continuous frames.
As a byproduct we give a partial answer to 
the question  raised by Ali, Antoine and Gazeau whether any continuous frame admits a corresponding discrete realization generated by sampling.
\end{abstract}




\section{Introduction}

In this paper we point out the relation between (continuous) frames and function spaces.
We illustrate that many function spaces can be described 
by 
continuous frames.
We further present a general method to derive atomic decompositions and Banach frames
for spaces with such a continuous frame description. Our results unify the theory of coorbit spaces
associated to integrable
group representations developed by Feichtinger and Gr\"ochenig \cite{FG1,FG2,FG3,gro} and its recent
generalizations \cite{DST,DST1,hr2}.

The concept of {\it discrete frames} in Hilbert spaces has been introduced by Duffin and
Schaeffer \cite{DS} and popularized greatly by Daubechies and her coauthors \cite{D2,DG}. 
A discrete frame is a countable family of elements in a separable Hilbert space
which allows stable not necessarily
unique (redundant) decompositions of arbitrary elements into 
expansions of the frame elements.
Later, motivated by the theory of coherent states, this concept was generalized by Antoine {\it et al.}
to families indexed by some locally compact space
endowed with a Radon measure. Their approach leads to the notion of continuous frames \cite{SAG,SAG1,GH,K}.
Prominent examples are connected to the continuous wavelet transform \cite{SAG,GM1}
and the short time Fourier transform \cite{Gr}.
In particular, {\it square integrable}
representations of groups generate continuous frames by
acting on a fixed mother atom. In mathematical physics, 
these frames are referred to as {\it coherent states} \cite{SAG,GM}.
Such decompositions into continuous superpositions of frame elements (atoms)
simplify the analysis of functions provided the atoms are suitably chosen with respect to the
problem under consideration.
For example, it is known that describing functions as continuous superposition of {\it wavelets} simplifies
the treatment of Calder\'on-Zygmund operators \cite{FJ2}, 
while  Gabor decompositions
{\it quasi-diagonalize} certain classes of pseudodifferential operators \cite{Gr}. 

Clearly, the concept of frame aims at stable decompositions in Hilbert spaces. However, 
in order to have a more complete and maybe finer characterization of reproducible vectors, 
one might ask whether these decompositions are also valid in certain Banach spaces. As a result of
this paper one may in fact 
associate to a continuous frame suitable Banach spaces, called {\it coorbit spaces}, provided the frame satisfies a certain integrability condition.
In these coorbit spaces, we have indeed stable decompositions. This might seem nearly 
a triviality at first glance because
the Banach space will be constructed in a way such that this is true. 
However, it turns out that in concrete examples
these Banach spaces are well-known classical function spaces, 
like
homogeneous or inhomogeneous
Besov or Triebel-Lizorkin spaces or 
the modulation spaces. In particular, these classes include the Sobolev spaces.

Usually in applications one 
prefers a discrete framework. So efforts have been done
to find methods to {\it discretize} classical continuous frames for
use in applications like signal processing, numerical solution of PDE, simulation,
and modeling, see for example \cite{SAG,DFR}.
In particular, the discrete wavelet transform and Gabor frames are prominent examples 
and have been proven to be a very successful tool for certain applications. Since the problem of discretization 
is so important it would be nice
to have a general method for this purpose. 
Indeed, Ali, Antoine and Gazeau 
asked for conditions which ensure that a certain sampling
of a continuous frame $\{\psi_x\}_{x\in X}$ yields a discrete frame $\{\psi_{x_i}\}_{i\in I}$
\cite[p.45]{SAG}.
As a byproduct of our results we give a partial answer to this question. In case the continuous
frame is generated by an integrable unitary representation of some group this is already covered
by results of Feichtinger and Gr\"ochenig \cite{FG1,FG2,FG3,gro}. 
Here, not only discrete frames for the corresponding
Hilbert space are constructed but at the same time Banach frames and atomic decompositions
for the associated coorbit spaces. By this general
theory it has been possible to unify atomic decompositions 
for important Banach spaces, like
homogeneous Besov-Triebel-Lizorkin spaces
\cite{FJ,Tr1,Tr2}, modulation spaces \cite{F4,Gr} and Bergman spaces \cite{FG1}.
A contribution by Antoine {\it et al.} \cite{AA}
describes another method to discretize continuous frames generated by {\it square}-integrable representations
of semidirect product groups $V\rtimes S$ where $V$ is a vector space and $S\subset GL(V)$ is a semisimple
connected Lie group.

Recently, it has been recognized that there exist several {\it continuous frames} with relevant
applications, which do not arise from some square integrable representation of
a group in a strict sense. So generalized concepts of coherent states have been introduced, 
where the continuous
frame is indexed by a homogeneous space $\G/H$ \cite{SAG}.
Important examples can be described in this setting, 
such as continuous wavelet and  Gabor frames on spheres \cite{AV,T3} and
continuous mixed Gabor / wavelet frames, i.e., continuous frames 
associated to the affine Weyl Heisenberg group
\cite{CF,HN,T1,T2,HL2,ff,f1}.
As a matter of fact, the theory of Feichtinger and Gr\"ochenig is no longer applicable in this setting.
So efforts to adapt their original approach to homogeneous spaces have been done recently by Dahlke
{\it et al.} \cite{DST,DST1}. For instance, they were able to define {\it modulation spaces}
on spheres as coorbit spaces. However, since they assume the continuous frame to be tight
their approach cannot currently cover most of the other cited examples.
Moreover, there are other examples of continuous frames which are not
indexed neither by groups nor homogeneous spaces. 
For instance, Rauhut \cite{hr1,hr2} constructed
continuous frames whose elements are invariant under the action of some symmetry groups. Hereby,
the corresponding index set is a space of orbits of the group under some compact automorphism group.
These frames could be used to describe subspaces of classical coorbit spaces consisting of elements which
are invariant under some symmetry group. 
Examples include subspaces of homogeneous Besov and Triebel-Lizorkin spaces and
modulation spaces of radially symmetric distributions. 
In particular, Banach frames and atomic decompositions of these spaces could be derived, where
all frame elements (atoms) are itself radial.

In this paper we introduce an abstract and generalized version of the Feichtinger and Gr\"ochenig approach
which unifies all earlier contributions \cite{DST, DST1,FG1,FG2,FG3,gro,hr2}.
The terminology of {\it coorbit spaces} for Banach spaces defined as retract of suitable solid Banach spaces 
by general transformations has firstly been introduced by Peetre \cite[p. 200]{P}. 
On one hand, our formulation is very much in the spirit of this initial concept, 
on the other hand it preserves the concrete applicability of the Feichtinger and Gr\"ochenig approach. 
In fact, we expect that this setting allows the characterization
of some other interesting function spaces, for example $\alpha$-modulation
spaces \cite{G,ff,HL2,HN}, Besov-Triebel-Lizorkin and modulation spaces on manifolds, e.g.,~on spheres.
The application of the present theory to these cases will be discussed elsewhere in successive contributions.

As a starting point we assume to have a general {\it continuous frame} 
$\mathcal{F}=\{\psi_x\}_{x \in X}$ indexed by
some locally compact space $X$.
We show that if the Gramian kernel $R(x,y)=\langle \psi_x, \widetilde \psi_y \rangle$ of the continuous
frame with respect to its {\it canonical dual} belongs to a certain Banach algebra of
integrable kernels on $X \times X$ then one can associate two classes of corresponding Banach spaces,
which we call (generalized) coorbit spaces.
We show that under certain {\it localization conditions} these two classes coincide.

As already announced we will describe a general method to sample a discrete set $(x_i)_{i\in I}\subset X$
such that $\mathcal{F}_d=\{\psi_{x_i}\}_{i \in I}$
is in fact a Banach frame or an atomic decomposition for the (generalized) coorbit spaces.
This is our main result and a further insight into the relations between continuous frames
and corresponding discrete frames in the description of Banach spaces.
It is recently recognized that {\it good} discrete frames for application and numerical implementation 
should have nice localization properties \cite{fg,Gr,gro1}. Indeed, we
are able to show that starting with a localized continuous frame 
the discrete frame arising from our sampling method is indeed localized in a suitable
sense.

The paper is organized as follows. Section 2 introduces some basic facts about continuous frames and our specific
assumptions. Two classes of coorbit spaces associated to a continuous frame and its canonical dual are described in
Section 3. Localization of continuous frames  and how localization ensures the coincidence of the two classes of
coorbit spaces are presented in Section 4. Section 5 is devoted to the discretization machinery. In particular,
we introduce the additional conditions on the continuous frame under which we can sample a discrete frame.
We conclude
the section by showing that the frames are in fact Banach frames for the original coorbit spaces.
In Section 6 we show that the discretization method preserves localization properties.
Finally, Section 7 is devoted to examples.

\section{Preliminaries}

Assume $\H$ to be a separable Hilbert space and $X$ a  locally compact Hausdorff space endowed with
a positive 
Radon measure $\mu$ with $\supp \mu = X$. For technical reasons we assume (without loss of generality)
that $X$ is $\sigma$-compact. In the following we denote  
generic constants, whose exact value is not important for a qualitative analysis, 
by $0<C,C',C'',C_1,C_2< \infty$.

A family $\F= \{\psi_x\}_{x\in X}$ of vectors in $\H$ is called a {\it continuous frame}
if there exist constants $0<C_1,C_2 < \infty$ such that
\beq\label{frame_ineq}
C_1 \|f\|^2 \leq \int_X |\langle f, \psi_x \rangle|^2 d\mu(x) \leq C_2 \|f\|^2\qquad \mbox{ for all } f\in \H.
\eeq
If $C_1 = C_2$ then the frame is called tight.
For the sake of simplicity we assume that the mapping $x \mapsto \psi_x$
is weakly continuous. Note that if $X$ is a countable set and $\mu$ the counting measure then 
we obtain the usual definition of a (discrete) frame.


Associated to $\F$ is the {\it frame operator} 
$S=S_{\cF}$ defined in weak sense by
\[
S :\H \to \H,\qquad S f \,:=\, \int_X \langle f, \psi_x \rangle \psi_x\, d\mu(x).
\]
From the {stability condition} (\ref{frame_ineq}) it follows that $S$ is a bounded, positive,
and boundedly invertible operator. If $\F$ is tight then $S$ is a multiple of the identity.
Furthermore, it also follows 
from (\ref{frame_ineq}) that
the set $\F$ is total in $\H$, i.e., $\F^\bot = \{0\}$, see also \cite{SAG1}.
We define the following two transforms associated to $\mathcal{F}$,
\begin{align}
V : \H \to L^2(X,\mu),\quad V f(x) \,:=&\, \langle f, \psi_x\rangle, \notag\\
W : \H \to L^2(X,\mu),\quad W f(x) \,:=&\, \langle f, S^{-1} \psi_x \rangle \,=\, V (S^{-1} f)(x).\notag
\end{align}
Their adjoint operators are given weakly by
\begin{align}
V^* : L^2(X,\mu) \to \H,\quad V^* F \,:=&\, \int_X F(y) \psi_y d\mu(y), \notag\\
W^* : L^2(X,\mu) \to \H,\quad W^* F \,:=&\, \int_X F(y) S^{-1} \psi_y d\mu(y). \notag
\end{align}
It holds $S = V^* V$, $S^{-1}=W^* W$, and $\id = V^* W = W^* V$.
In fact, since $S$ is invertible and self-adjoint we have
\beq\label{W_invert}
f \,=\, S S^{-1} f \,=\, \int_X \langle S^{-1} f, \psi_y \rangle \psi_x d\mu(y)
\, = \, \int_X W f(y) \psi_y d\mu(y)
\eeq
in weak sense. Of course, this is an inversion formula for $W$. Replacing $f$ by $Sf$ yields
an inversion formula
for $V$, i.e.,
$
f = \int_X V f(y) S^{-1} \psi_y d\mu(x).
$
Forming the scalar product with $S^{-1} \psi_x$ in (\ref{W_invert}), resp. with $\psi_x$ in the inversion formula for $V$,
 yields
\[
W f(x) \,=\, \int_X Wf(y) \langle \psi_y, S^{-1} \psi_x \rangle d\mu(y)
\quad \mbox{and} \quad 
V f(x) \,=\, \int_X Vf(y) \langle \psi_y, S^{-1} \psi_x \rangle d\mu(y).
\]
Hence, it make sense to define the kernel
\begin{equation}\label{def_kernelR}
R(x,y) \,:=\, R_\F(x,y) \,:=\,  \langle \psi_y, S^{-1} \psi_x \rangle.
\end{equation}
Denoting the application of a kernel $K$ to a function $F$ on $X$ by
\begin{equation}\label{def_kernel_app}
K(F)(x) \,:=\, \int_X F(y) K(x,y) d\mu(y)
\end{equation}
we have $R(W f) = W f$ and $R(Vf) = Vf$ for all $f \in \H$.
Since $S$ is self-adjoint,
it holds $R(x,y) = \ol{R(y,x)}$. This means that $R$ is self-adjoint as an operator on $L^2(X,\mu)$.
Furthermore, the mapping $F \mapsto R(F)$ is an orthogonal projection from
$L^2(X,\mu)$ onto the image of $W$ (which equals the image of $V$).

If $\tilde{\cF}=\{\tilde{\psi}_x\}_{x\in X}$ is another frame that satisfies
\begin{equation}\label{dual2}
        f = \int_X \langle f, \psi_x\rangle \tilde \psi_x d\mu(x), \quad \text{for all } f \in \cH
\end{equation}
then $\tilde{\cF}$ is called a {\it dual frame}. In particular, $S^{-1} \F = \{S^{-1} \psi_x\}_{x\in X}$
is a dual frame, called the {\it canonical dual frame}. Since in general $\text{ker}(V^*) \neq \{0\}$ there
may exist several dual frames.

We assume in the following that $\|\psi_x\| \leq C$ for all $x \in X$. This implies by the Cauchy Schwarz inequality
$|V f(x)| \leq C \|f\|$ and $|W f(x)| \leq C\on S^{-1}\on\,\|f\|$ for all $x\in X$ and, together with the weak continuity
assumption, we conclude $Vf, Wf \in C^b(X)$ for all $ f\in \H$, where $C^b(X)$ denotes the bounded continuous functions
on $X$. 

In the sequel we denote by $\cB(Y)$ the bounded linear operators from a Banach space into itself and by
$\on \cdot | Y\on$ its norm.

\section{Coorbit Spaces}

Inspired by the pioneering work of Feichtinger and Gr\"ochenig {\it et al.} \cite{FG1,FG2,FG3,gro,gro1,fg} and
their recent generalizations \cite{DST,DST1,hr1,hr2}, we show in this section how classes of abstract Banach spaces
called (generalized) {\it coorbit spaces} can be associated to a given continuous frame.
Such Banach spaces will essentially describe vectors (or functionals) $f$ such that the corresponding
transforms $V f$ or $W f$ belongs to a fixed parameter space $Y$ of functions on $X$. In concrete examples
the coorbit spaces are certain function spaces.

In order to detail this idea, we need some preparation. We define the 
Banach algebra of kernels 
\[
\A_1 \,:= \, \{K:X\times X \to \C,~K \mbox{ measurable }, \|K|\A_1\| < \infty\}
\]
where
\[
\|K|\A_1\| \,:=\,
\max\left\{ \esssup_{x\in X} \int_{X} |K(x,y)| d\mu(y), \esssup_{y\in X} \int_X |K(x,y)| d\mu(x)
\right\}
\]
is its norm. The multiplication in $\A_1$ is given by
\begin{equation}\label{A_mult}
K_1 \circ K_2 (x,y) \,=\, \int_X K_1(x,z) K_2(z,y) d\mu(z).
\end{equation}
Identifying a kernel with an operator acting on suitable functions $F$ on $X$ by (\ref{def_kernel_app})
we clearly have $K_1(K_2(F)) = (K_1 \circ K_2)(F)$.

In the sequel we make the basic assumption that the kernel
$R$ defined in (\ref{def_kernelR}) is contained in $\A_1$.


We will also need suitable weighted subalgebras of $\A_1$. We call a {\it weight function} $m:X \times X \to \R$
{\it admissible} if $m$ is continuous,
\begin{align}
\label{m_submult}
1\,\leq\, m(x,y) \,\leq&\, m(x,z) m(z,y)\quad \mbox{ for all }x,y,z \in X,\\
\label{m_symmetry}
m(x,y) \,=&\, m(y,x) \quad \mbox{ for all } x,y \in X, \\
\label{m_diag_bound}
\mbox{and }\quad  m(x,x) \,\leq&\, C \,<\, \infty \mbox{ for all } x\in X.
\end{align}
For an admissible weight $m$ we define the Banach algebra
\[
\A_m := \{K: X\times X \to \C,~ Km \in \A_1\}
\]
endowed with the natural norm
\[
\|K|\A_m\| \,:=\, \|Km|\A_1\|.
\]
Property (\ref{m_submult}) ensures that $\A_m$ is in fact an algebra with the multiplication (\ref{A_mult}). Moreover,
the symmetry property (\ref{m_symmetry}) ensures that $\A_m$ is a Banach-$*$-algebra with the involution
$K^*(x,y) = \ol{K(y,x)}$. Interpreting $K$ as an operator on $L^2(X,\mu)$ its adjoint is in fact given by the kernel $K^*$.

In order to define our {\it coorbit spaces} associated to the continuous frame, we will make use
of a function space $Y$ that satisfies the following properties.
\begin{itemize}
\item[(Y1)] $(Y,\|\cdot|Y\|)$ is a non-trivial Banach space of functions on $X$ that is continuously embedded into $L^1_{loc}(X,\mu)$,
and that satisfies the {\it solidity condition}, i.e., if $F \in L^1_{loc}(X,\mu), G\in Y$, such that
$|F(x)| \leq |G(x)|$ a.e. then $F \in Y$ and $\|F|Y\| \leq \|G|Y\|$.
\item[(Y2)] There exists an admissible weight function $m$ such that $\A_m (Y) \subset Y$ and
\begin{equation}\label{KF_ineq}
\|K(F)|Y\| \,\leq\, \|K|\A_m\|\,\|F|Y\| \quad \mbox{ for all } K\in \A_m, F \in Y.
\end{equation}
\end{itemize}
By Schur's test (sometimes also referred to as generalized Young inequality)
\cite[Theorem 6.18]{Folland2} all $L^p(X,\mu)$ spaces, $1 \leq p \leq \infty$, are
examples for such $Y$ spaces (with
trivial weight $m=1$).  Moreover, if $w$ is a continuous weight function on $X$,
we define $L^p_w:=L^p_w(X,\mu) =\{F, Fw \in L^p(X,\mu)\}$ with norm $\|F|L^p_w\|:=\|Fw|L^p\|$ and denote
\begin{equation}\label{mw_associate}
m(x,y) \,:=\, \max\left\{\frac{w(x)}{w(y)},\frac{w(y)}{w(x)}\right\}.
\end{equation}
Then $m$ is admissible and $\A_m (L^p_w) \subset L^p_w$ again by Schur's test.

In the following we only admit $Y$ with properties (Y1) and (Y2) and such that $R$ defined in 
(\ref{def_kernelR}) is contained in $\A_m$, where
$m$ is the weight function associated to $Y$.

The next step is to derive a reservoir 
to embed our Banach spaces in. 
To this end take a fixed point $z \in X$ and define a weight function on $X$ by
\beq\label{v_def}
v(x)\,:=\,v_z(x)\,:=\,m(x,z).
\eeq
By the properties of $m$, the choice of another point $z'$ yields an {\it equivalent weight}, i.e.,
$v_{z'}(x) = m(x,z') \leq m(x,z) m(z,z') = m(z,z') v_z(x)$. Exchanging the roles of $z$ and $z'$ gives
a reversed inequality. Now, we define the spaces
\begin{align}
\H_v^1 \,:=&\, \{f \in \H, V f \in L^1_v\},
\qquad \K_v^1 \,:=\,\, \{f \in \H, W f \in L^1_v\} \notag
\end{align}
with natural norms
\begin{align}\label{def_normH}
\|f|\H^1_v\| \,:=\, \|V f| L^1_v\|,
\qquad
\|f|\K^1_v\| \,:=\, \|W f| L^1_v\|.
\end{align}
Since $\{\psi_x\}_{x\in X}$ is total in $\H$, the expressions in (\ref{def_normH}) 
indeed define norms, not only seminorms.
The operator $S$ is obviously an isometric isomorphism between $\H^1_v$ and $\K^1_v$.

\begin{proposition} The spaces $(\H^1_v,\|\cdot|\H^1_v\|)$ and $(\K^1_v,\|\cdot|\K^1_v\|)$ are Banach spaces.
\end{proposition}
\begin{Proof} Suppose that $(f_n)_{n \in \N} \subset \H$ is a Cauchy sequence in $\H^1_v$. This means that $(F_n)=(V f_n)$
is a Cauchy sequence in $L^1_v$ and by completeness of $L^1_v$ it holds $F_n \to F \in L^1_v$.
Furthermore, it holds $R(F_n) = F_n$ by the reproducing formula on the image of $\H$ under $V$.
This implies that $R(F) = F$. Since 
$|R(x,y)| \leq C^2\on S^{-1}|\H \on$
for all $x,y \in X$ and
$v(x) \geq 1$ 
it holds
\[
|R(F)(x)| \leq \int_X |F(y)| |R(x,y)| d\mu(y)
\,\leq\, \int_X |F(y)| v(y) |R(x,y)| d\mu(y) \leq C^2\on S^{-1}\on\,\|F|L^1_v\|
\]
implying $F=R(F) \in L^\infty$. By $L^\infty \cap L^1_v \subset L^2$ we have
$F = R(F) \in L^2(X,\mu)$. Since the application of $R$ is the orthogonal projection from $L^2$ onto
the image of $V$, there exists $f \in \H$ such that $F = V f$. Moreover, $Vf \in L^1_v$ means $f \in \H^1_v$
and $f_n \to f$ in $\H^1_v$.
The same arguments show that also $\K^1_v$ is a Banach space.
\end{Proof}

Since $R$ is assumed to be in $\A_m$ we obtain
\begin{align}
\|\psi_y|\K^1_v\| \,=&\, \int_X |W (\psi_y)(x)| v(x) d\mu(x)
\,=\, \int_X |R(x,y)| m(x,z) d\mu(x) \notag\\
\label{Knorm_psix}
\,\leq& \, m(y,z) \int_X |R(x,y)| m(x,y) d\mu(x) \leq
v(y) \|R|\A_m\|.
\end{align}
and similarly
\beq\label{norm_psix}
\|S^{-1}\psi_y|\H^1_v\| \,=\, \int_X |V (S^{-1}\psi_y)(x)| v(x) d\mu(x) \,\leq\, v(y) \|R|\A_m\|.
\eeq
Hence, $\psi_y \in \K^1_v$ and $S^{-1} \psi_y \in \H^1_v$ for all $y \in X$. Denote by $(\H^1_v)^\urcorner$
the space of all continuous conjugate-linear functionals on $\H^1_v$ (the anti-dual) and similarly define
$(\K^1_v)^\urcorner$. 
We extend the bracket on $\H$ to $(\H^1_v)^\urcorner \times \H^1_v$
by $\langle f,g\rangle = f(g)$ for $f \in (\H^1_v)^\urcorner$, $g \in \H^1_v$ and analogously for
$(\K^1_v)^\urcorner \times \K^1_v$. Taking the anti-dual instead of the dual yields the same
calculation rules for the bracket as in the Hilbert space setting.
Since $\spann\{\psi_x, x\in X\}$ and $\spann\{S^{-1} \psi_x, x\in X\}$ are dense in $\H$ the spaces
$\H^1_v$ and $\K^1_v$ are dense in $\H$ and
$\H$ is weak-$*$ dense in $(\H^1_v)^\urcorner$ and in $(\K^1_v)^\urcorner$.

Since $\psi_x \in \K^1_v$ we may extend the transform $V$ to $(\K^1_v)^\urcorner$ by
\[
V f(x) \,=\, \langle f, \psi_x\rangle \,=\, f(\psi_x), \quad f\in (\K^1_v)^\urcorner.
\]
By the same argument, the transform $W$ extends to $(\H^1_v)^\urcorner$,
\[
W f(x) \,=\, \langle f, S^{-1} \psi_x\rangle \,=\, f(S^{-1}\psi_x), \quad f \in (\H^1_v)^\urcorner.
\]
We may also extend the operator $S$ to an isometric
isomorphism between $(\K^1_v)^\urcorner$ and
$(\H^1_v)\urcorner$ by $\langle Sf, g\rangle = \langle f, Sg\rangle$ for $f \in (\K^1_v)^\urcorner$
and $g \in \H^1_v$ (recall that $Sg \in \K^1_v$).

Next, we need to show that $\spann\{\psi_x, x\in X\}$ 
and $\spann\{S^{-1}\psi_x, x\in X\}$ are dense in $\K^1_v$ and
$\H^1_v$, respectively.

\begin{lemma}\label{lem_Hdual} The expression $\|Vf|L^\infty_{1/v}\|$ is an equivalent norm on 
$(\K^1_v)^\urcorner$ and
$\|Wf|L^\infty_{1/v}\|$
is an equivalent norm on $(\H^1_v)^\urcorner$.
\end{lemma}
\begin{Proof} Observe that for $f \in (\K^1_v)\urcorner$ it holds by (\ref{Knorm_psix})
\[
|V f(x)| \,=\, |\langle f, \psi_x\rangle| \,\leq\, \|f|(\K^1_v)^\urcorner\|\,\|\psi_x|\K^1_v\|
\leq \|f|(\K^1_v)^\urcorner\|\,\|R|\A_m\| v(x).
\]
For the converse we use that $V^*$ is an isometric isomorphism from $R(L^1_v)$ to $\K^1_v$, i.e.,
\begin{align}
\|f|(\K^1_v)^\urcorner\| \,=&\, \sup_{\|h|\K^1_v\|=1} |\langle f,h\rangle|
\,=\, \sup_{H \in R(L^1_v),\|H|L^1_v\|\leq 1} |\langle f, V^* H\rangle| \notag\\
\,\leq&\, \sup_{H \in L^1_v, \|H|L^1_v\|\leq 1} |\langle Vf, H\rangle|
\,=\, \|Vf | L^\infty_{1/v}\|.\notag
\end{align}
The assertion for $\H^1_v$ is shown in the same way.
\end{Proof}

\begin{Theorem}\label{thm_frame_dense} \begin{itemize}\item[(a)] $\K^1_v$ is characterized by the vector space of all uniform unconditional expansions of the form
\begin{equation}\label{f_repr}
f \,=\, \sum_{i\in I} c_i \psi_{x_i}
\end{equation}
where $(x_i)_{i\in I}$ is an arbitrary countable subset of $X$ and
\[
\sum_{i \in I} |c_i| v(x_i) < \infty.
\]
The expression
\[
\|f\|':=\inf \sum_{i\in I} |c_i| v(x_i)
\]
where the infimum is taken over all representation (\ref{f_repr}) of $f$, is an equivalent norm on $\K^1_v$. In particular,
$\spann\{\psi_x, x\in X\}$ is dense in $\K^1_v$.
\item[(b)] $\H^1_v$ is characterized analogously by replacing $\psi_{x_i}$ by $S^{-1} \psi_{x_i}$ in (\ref{f_repr})
\end{itemize}
\end{Theorem}
\begin{Proof} The proof is completely analogous to the proof of Theorem 12.1.8 in \cite{Gr} and hence omitted, see also \cite{Bonsall,hr}.
We remark only that Lemma \ref{lem_Hdual} is used.
\end{Proof}

It follows that $\H^1_v$ and $\K^1_v$ have a certain minimality property.

\begin{corollary}\label{cor_min} Suppose that $(B,\|\cdot|B\|)$ is a Banach space that contains all frame elements $\psi_x, x\in X,$ and assume
that $\|\psi_x|B\| \leq C v(x)$ for some constant $C$. Then $K^1_v$ is continuously embedded into $B$. The same holds
replacing $\psi_x$ by $S^{-1} \psi_x$ and $\K^1_v$ by $\H^1_v$.
\end{corollary}
\begin{Proof} 
%
%
Suppose that $f=\sum_{i\in I} c_i \psi_{x_i}$ with $\sum_{i\in I} |c_i| v(x_i) < \infty$, i.e., $f \in \K^1_v$ by
Theorem \ref{thm_frame_dense}. Then
\[
\|f|B\| \leq \sum_{i\in I} |c_i| \|\psi_{x_i}|B\| \leq C\sum_{i\in I} |c_i| v(x_i) < \infty.
\]
This implies $f \in B$. Taking the infimum over all possible expansions of $f$ yields $\|f|B\| \leq C \|f|\K^1_v\|$ and
the embedding is continuous.
\end{Proof}

Let us now give a more precise statement about the weak-$*$ density of $\H$ in $(\H^1_v)^\urcorner$
and in $(\K^1_v)^\urcorner$, respectively.


\begin{lemma}\label{lem_wstar} Let $f \in (\H^1_v)^\urcorner$. Then there exists a sequence 
$(f_n)_{n \in \N} \subset \H$ with $\|f_n|(\H^1_v)^\urcorner\| \leq C \|f|(\H^1_v)\urcorner\|$
such that $f_n$ is weak-$*$ convergent to $f$. The same holds with $(\H^1_v)^\urcorner$ replaced by 
$(\K^1_v)^\urcorner$.
\end{lemma}
\begin{Proof} Since $X$ is $\sigma$-compact there exists a sequence of compact subsets $(U_n)_{\in \N}$
such that $U_{n} \subset U_{n+1}$ and $X = \bigcup_{n = 1}^\infty U_n$. Denote by $\chi_{U_n}$ the 
characteristic function of $U_n$, let 
$F_n = \chi_{U_n} Wf \in L^2(X,\mu)$ and
set $f_n = V^* F_n = \int_{U_n} Wf(y) \psi_y d\mu(y) \in \H$, $n \in \N$. 
It is straightforward to show that the sequence $(f_n)_{n \in \N}$ has the desired properties.
\end{Proof}

\begin{lemma}\label{lem_HK} \begin{itemize}
\item[(a)] For $f\in (\K^1_v)^\urcorner$ it holds $Vf \in L^\infty_{1/v}(X)$ and for
$f \in (\H^1_v)^\urcorner$ it holds $Wf \in L^\infty_{1/v}(X)$. The mappings
$V : (\K^1_v)^\urcorner \to L^\infty_{1/v}(X)$ and $W: (\H^1_v)^\urcorner \to L^\infty_{1/v}(X)$
are injective.
%
%
%
%
\item[(b)] A bounded net $(f_\alpha)_{\alpha \in I}$ in $(\K^1_v)^\urcorner$ (resp. in $(\H^1_v)^\urcorner)$
is weak-$*$ convergent
to an element $f \in (\K^1_v)^\urcorner$ (resp. $f \in (\H^1_v)^\urcorner$)
if and only if $Vf_\alpha$ (resp. $W f_\alpha$) converges
pointwise to $Vf$ (resp. $Wf$).
\item[(c)] The reproducing formula extends to
$(\K^1_v)^\urcorner$ and $(\H^1_v)^\urcorner$, i.e.,
\begin{align}
V f \,=&\, R(Vf) \qquad \mbox{ for all } f\in (\K^1_v)^\urcorner,\\
W f \,=&\, R(Wf) \qquad \mbox{ for all } f\in (\H^1_v)^\urcorner.
\end{align}
\item[(d)] Conversely, if $F \in L^\infty_{1/v}(X)$ satisfies the reproducing formula $F = R(F)$
then there exist $f \in (\K^1_v)\urcorner$ and $g \in (\H_v^1)^\urcorner$ such that $F = V f = W g$.
\end{itemize}
\end{lemma}
\begin{Proof} (a) The assertion follows from Lemma \ref{lem_Hdual}.

(b) The assertion follows from the density of $\spann \{\psi_x, x\in X\}$ (resp. $\spann\{S^{-1}\psi_x, x\in X\}$)
in $\K^1_v$ (resp. $\H^1_v$) and by definition of $Vf$ (resp. $Wf$).

(c) Suppose that $f \in (\K^1_v)^\urcorner$. Then by Lemma~\ref{lem_wstar} 
there exists a sequence $(f_n)_{n \in \N}\subset \H$, which is
weak-$*$ convergent to $f$ and norm bounded in $(\K^1_v)^\urcorner$. This implies
the pointwise convergence of $V f_n$ to $Vf$. Since $f_n \in \H$ 
the reproducing formula
holds for $V f_n$. Furthermore, we have 
$|V f_n(y)| \leq C \sup_{n \in \N} \|f_n|(\K^1_v)^\urcorner\| v(y) \leq C' v(y)$
and $y \mapsto v(y) R(x,y)$ is
integrable for any $x \in X$. Thus,
Lebesgue's dominated convergence theorem yields
\[
V f(x) = \lim_{n \to \infty} Vf_n(x) = \lim_{n \to \infty} \int_X R(x,y) V f_n(y) d\mu(y)
= \int_X R(x,y) Vf(y) d\mu(y) = R(Vf)(x).
\]
Analogously we obtain the reproducing formula for $W f$.

(d) A direct computation shows that the adjoint mappings of $V:\H^1_v \to L^1_v$ and $W: \K^1_v\to L^1_{v}$
are given weakly by
\begin{align}
V^*: L^\infty_{1/v} \to (\H^1_v)^\urcorner,\quad&
V^* F \,=\, \int_X F(x) \psi_x d\mu(x) \quad\mbox{ for } F \in L^\infty_{1/v},\notag\\
W^*: L^\infty_{1/v} \to (\K^1_v)^\urcorner,\quad&
W^* F \,=\, \int_X F(x) S^{-1} \psi_x d\mu(x) \quad\mbox{ for } F \in L^\infty_{1/v}.\notag
\end{align}
It holds
\[
W (V^*F)(y) \,=\, \int_X F(x) \langle\psi_x, S^{-1} \psi_y\rangle d\mu(x)
\,=\, R(F)(y)
\]
and similarly $V(W^* F) = R(F)$. Hence, if $F=R(F)$ then $F= V f = Wg$, where
$f=W^*F \in (\K^1_v)^\urcorner$ and $g = V^* F \in (\H^1_v)^\urcorner$.
\end{Proof}

Now we are ready to define the coorbit spaces.

\begin{definition}\label{def_coorbit} The coorbits of $Y$ with respect to the frame $\F=\{\psi_x\}_{x\in X}$ 
are defined as
%
\begin{align}
\Co Y  \,:=\, \Co(\F,Y) \,:=&\, \{f \in (\K^1_v)^\urcorner, V f \in Y\},\notag\\
\TCo Y \,:=\, \Co(S^{-1} \F,Y) \,:=&\, \{f \in (\H^1_v)^\urcorner, W f \in Y\}\notag
\end{align}
with natural norms
\begin{align}
\|f | \Co Y\| \,:=\, \|V f| Y\|,\qquad \|f | \TCo Y\|\,:=\, \|W f| Y\|.\notag
\end{align}
\end{definition}
Obviously, the operator $S$ is an isometric isomorphism
between $\Co Y$ and $\TCo Y$.

\begin{proposition}\label{Prop_proj} Suppose that $R(Y) \subset L^\infty_{1/v}$.
Then the following statements hold.
\begin{itemize}\itemsep-1pt
\item[(a)] The spaces $(\Co Y,\|\cdot|\Co Y\|)$ and $(\TCo Y,\|\cdot|\TCo Y\|)$ are Banach spaces.
\item[(b)] A function $F \in Y$ is of the form $Vf$ (resp. $W f$) for some $f \in \Co Y$ (resp. $f \in \TCo Y$)
if and only if $F = R(F)$.
\item[(c)] The map $V: \Co Y \to Y$ (resp. $W : \TCo Y \to Y$) establishes an isometric isomorphism between
$\Co Y$ (resp. $\TCo Y$)
and the closed subspace $R(Y)$ of $Y$.
\end{itemize}
\end{proposition}
\begin{remark} The condition $R(Y) \subset L^\infty_{1/v}$ might seem a bit strange at
first glance. However, we will show later in Corollary \ref{cor_RY} that, under the assumption
we will need to develop the discretization method in Section 5, this property holds true automatically.

Of course, if one wants to be sure that the spaces $\Co Y$ and $\TCo Y$ are Banach spaces in any case,
one may take the norm or weak-$*$ completion of the spaces in Definition \ref{def_coorbit}.
\end{remark}
\begin{Proof} Let us first prove (b). If $f \in \Co Y$ then by definition
$f \in (\K^1_v)^\urcorner$ and hence by Lemma \ref{lem_HK}(c) the reproducing formula holds.
Conversely, if $F \in Y$ satisfies $R(F) = F$ then $F \in L^\infty_{1/v}$ by the assumption
$R(Y) \subset L^\infty_{1/v}$. Lemma \ref{lem_HK}(d) implies that there exists 
$f \in (\K^1_v)^\urcorner$ such that $F = Vf$. Since $F \in Y$ we have $f\in \Co Y$.
Analogous arguments apply to $\TCo Y$.

(a) Suppose that $(f_n) \subset \Co Y$ is a Cauchy sequence implying that $F_n = V f_n$ is a Cauchy
sequence in $Y$ converging to an element $F \in Y$ by completeness of $Y$. By (b) it holds
$F_n = R(F_n)$ and since $R$ acts continuously on $Y$ we conclude $R(F) = F$.
Again by (b) there exists hence an $f \in \Co Y$ such that $F = V f$.
The analogous arguments apply to $\TCo Y$.

(c) The assertion follows from the injectivity of $V$ and $W$ (Lemma \ref{lem_HK}(a)) together
with (b).
\end{Proof}

\begin{corollary} 
\begin{itemize}\itemsep-1pt
\item[(a)] It holds $\Co L^\infty_{1/v} = (\K^1_v)^\urcorner$ and
$\TCo L^\infty_{1/v} = (\H^1_v)^\urcorner$.
\item[(b)] It holds $\Co L^2 = \TCo L^2 = \H$.
\item[(c)] Suppose $R(Y) \subset L^\infty_{1/v}$. Then the definition of the
coorbit spaces is independent of the weight function $m$,
resp. $v$, i.e., if $m_2$ is another weight with $m(x,y) \leq C m_2(x,y)$ and $v_2$ is the
corresponding weight function on $X$ then
\begin{align}
\Co Y \,=&\, \{f \in (\K^1_{v_2})^\urcorner, Vf \in Y\},\notag\\
\TCo Y \,=&\, \{f \in (\H^1_{v_2})^\urcorner, Wf \in Y\}.\notag
\end{align}
\end{itemize}
\end{corollary}
\begin{Proof} The statement (a) was already proved in Lemma \ref{lem_Hdual}, and (b), (c) are
shown as in \cite{FG1}.
\end{Proof}


\begin{remark}\label{rem_Y2} Analyzing the proofs of this section one might realize that the assumption ({\rm Y2})
on the function space $Y$ may be weakened. Actually we only needed that $R$ is contained in $\A_m$ and
that $R$ acts continuously on $Y$. 
So in order to define
the coorbit spaces corresponding to $Y$
it is enough that the subalgebra $\A:= \A_m \cap \cB(Y)$
is not trivial
and $R$ is contained in $\A$. Of course, if $Y$ is a weighted $L^p$ space and $m$ is the corresponding
weight (\ref{mw_associate}) then $\A$ coincides with $\A_m$. However, there are function spaces $Y$
for which $\A$ is a proper subalgebra of $\A_m$, for instance if $Y$ is a certain mixed norm space.
\end{remark}

\section{Localization of Frames}

It might seem strange at first glance that, for a given continuous frame, we have to deal with two classes of coorbit spaces.
So the question arises under which conditions it holds $\Co Y = \TCo Y$. 
Furthermore, it is interesting to  investigate the dependence of the coorbit spaces on the particular
frame chosen.

The main tool for these investigations will be the concept of {\it localization of frames}.
In particular, we will generalize the theory developed in \cite{gro1,gl,fg} of intrinsic localization of
discrete frames to continuous ones. We will show that if $\F$ is {\it intrinsically $\cA$-localized},
i.e., its {\it Gramian} kernel belongs to a suitable {\it spectral algebra} $\cA$, then also its
canonical dual $S^{-1} \F$ is intrinsically $\cA$-localized and $\Co Y = \TCo Y$. Moreover, we will
show that 
$\cA$-localization is an equivalence relation in the set of intrinsically $\cA$-localized frames and that
equivalent frames generate in fact equivalent coorbit spaces.

\subsection{$\cA$-localized Continuous Frames}


Let us first specify the algebras of kernels which are suitable to measure localization. 

\begin{definition}\label{admAlg}
 A  Banach-$*$-algebra $(\cA,\|\cdot|\cA\|) $ of kernels $K$ on $X \times X$ is called {\em admissible with respect to $(X,\mu)$}
if the following properties hold:
\begin{itemize}
\item[(A1)] $\cA$  is a continuously embedded into $\cB (L^2(X,\mu) ) $;
\item[(A2)] $\cA $ is {\it solid}, i.e., if $K$ is a measurable kernel, 
$|K| \leq |R|$ and $R \in \cA$ then also $K \in \cA$ and $\|K|\cA\| \leq \|R|\cA\|$.
\end{itemize}
The multiplication in the algebra is assumed again as in formula \eqref{A_mult}. Moreover we assume that the action of $K \in \cA$ on a function $F \in L^2(X,\mu)$ is given as in \eqref{def_kernel_app}. This again implies that $A_1(A_2(F)) = (A_1 \circ A_2) (F)$, for all $A_1,A_2 \in \cA$, and $F \in L^2(X,\mu)$. 
\end{definition}


In the following we assume $\cA$ to be admissible with respect to $(X,\mu)$. Of course, the algebras $\A_m$ from the previous section are admissible. 

Given two continuous frames $\cG=\{g_x\}_{x\in X},\cF=\{f_x\}_{x\in X}$ in $\cH$, their
{\it crossed Gramian kernel} is defined as
\[
G(\cG,\cF)(x,y):=\langle f_y,g_x \rangle.
\]

\begin{definition}\label{frame}
A frame $\cG$ for $\cH $ is called $\cA$-localized with respect to
a frame $\cF $ if $G(\cG,\cF) \in \cA $. In this case we write $\cG \sim_{\cA } \cF $.
If $\cG \sim_{\cA} \cG$, then $\cG $ is called $\cA $-self-localized or
intrinsically  $\cA $-localized.
\end{definition}
In the following $\tilde {\cF}$ always denotes a dual frame of the frame $\cF$.
Since $G(\cG,\cF)^* = G(\cF,\cG)$ and $\cA$ is assumed to be a $*$-algebra the relation
$\sim_{\cA}$ is symmetric.
One might ask whether $\sim_{\cA}$ is an equivalence relation. This is not true in general.
In fact, reflexivity holds only for intrinsically $\cA $-localized
frames. Transitivity is valid only in a modified version involving also dual frames as shown in
the following lemma.

\begin{lemma}\label{semieq}
Let $\cE =\{e_x\}_{x \in X}, \cF=\{f_x\}_{x \in X},\cG=\{g_x\}_{x \in X}$ be frames for $\mathcal{H}$.
\begin{itemize}
\item[(a)] If $\cE \sim_{\cA} \tilde{\cF}$ and $\cF \sim_{\cA} \cG $ then we have $\cE \sim_{\cA} \cG $.
\item[(b)] If $\cF \sim_{\cA} \cF $ and $\tilde{\cF} \sim_{\cA} \tilde{\cF}$ for a dual frame $\tilde{\cF}$
then it holds $\cF \sim_{\cA} \tilde{\cF} $.
\item[(c)] If both $\tilde{\cF} _1$ and $\tilde{\cF} _2 $ are dual frames of $\cF $,
 $\cF \sim_{\cA} \cF $ and $\tilde{\cF} _1 \sim_{\cA} \tilde{\cF} _2$ then 
$\cF \sim_{\cA} \tilde{\cF} _j$ for $j=1,2$.
\end{itemize}
\end{lemma}
\begin{Proof}
Let us show (a). The statements (b) and (c) 
are just direct consequences of (a).
Since
$e_y = \int_{X} \langle e_y, \tilde f_z \rangle f_z d\mu(z)$
we have
$$
\langle e_y, g_x\rangle \,=\, \int_{{X}} \langle e_y, \tilde f_z \rangle  \langle f_z,g_x\rangle d\mu(y)\quad\text{ for all } x,z \in X.
$$
This immediately implies that
$|G(\cE,\cG)| \leq   |G(\cF,\cG)| \circ |G(\cE,\tilde{\cF})|$
and one concludes by solidity of $\cA$.
\end{Proof}
\begin{remark}
If $\cF \sim_{\cA} \cF $ and likewise for the canonical dual frame $S^{-1}\cF \sim_{\cA} S^{-1} \cF$,
then Lemma \ref{semieq} (b) ensures that $R=R_\F = G(\F,S^{-1}\F) \in \cA$.
Of course, $\A=\A_m$ is of particular interest when considering coorbit spaces.
In Section \ref{sec_dual_loc} we will show that
$\cF \sim_{\cA} \cF $ implies 
$S^{-1}\cF \sim_{\cA} S^{-1} \cF$ 
under a certain assumption on the algebra $\A$.
\end{remark}

\subsection{Localization Conditions Ensure $\Co Y =\TCo Y$}

In the following we assume that $\Co Y$ and $\TCo Y$ are Banach spaces.
So in case we are in the pathological situation where the spaces in Definition
\ref{def_coorbit} are not complete we take their completion. As usual $m$ denotes the weight function
associated to $Y$.

\begin{proposition}\label{eqCo}
Suppose that $\cF$ is a frame for $\H$ with canonical dual $S^{-1}\F$.
If both $\F$ and $S^{-1}\F$ are intrinsically $\A_m$-localized then
$\Co Y=\TCo Y$ with equivalent norms. In particular, it holds $\cH^1_v=\cK^1_v$.
\end{proposition}
\begin{Proof}
For $f \in \cH$ we have
\begin{equation}
\label{WV}
        Wf\,=\, G(S^{-1}\F, S^{-1}\F) (Vf) \qquad \mbox{ and } \qquad Vf = G(\cF,\cF) (Wf).
\end{equation}
By (\ref{KF_ineq}) this implies $\cH^1_v=\cK^1_v$ and $(\cH^1_v)^\urcorner=(\cK^1_v)^\urcorner$.
Since formulae \eqref{WV} extend to $(\cH^1_v)^\urcorner$ and $\cA_m(Y) \subset Y$ 
(Lemma~\ref{lem_HK}(c))
we immediately obtain $\|Wf|Y\| \asymp \| V f|Y\|$.
\end{Proof}

In the following we will show that for an $\cA_m$-self-localized frame  the space $\Co Y$
can be characterized by using different $\cA_m$-self-localized duals and
that any other $\cA_m$-self-localized frame $\cG$ which is localized to $\cF$ generates 
the same spaces.

\begin{proposition}\label{eqCo1}
Assume that $\mathcal{F}=\{\psi_x\}_{x\in X}$, $\tilde{\mathcal{F}}=\{\tilde \psi_x\}_{x \in X}$ 
are two mutually dual
$\cA_m$-self-localized frames for $\mathcal{H}$. 
Then it holds $\Co(\tilde \cF,Y) = \TCo Y$
with equivalent norms. In particular, if the canonical dual of $\mathcal{F}$ is $\cA_m$-self-localized,
then $\Co(\tilde \cF,Y) = \TCo Y = \Co Y$.
\end{proposition}
\begin{Proof}
Since $\mathcal{F}$ and $\tilde{\mathcal{F}}$ are mutually dual $\cA_m$-self-localized frames
one obtains with Lemma \ref{semieq}(b) that $\mathcal{F} \sim_{\cA_m} \tilde \cF$ and
therefore $\tilde{\mathcal{F}} \subset \cH^1_v$. Denote 
$\tilde V f(x):=\langle f, \tilde \psi_x \rangle$.
By expanding the elements $S^{-1} \psi_x$ of the
canonical dual frame with respect to $\tilde \psi_x$ by using formula \eqref{dual2} and then $\tilde \psi_x$
with respect to $S^{-1} \psi_x$ we obtain
\[
        W f \,=\, R(\tilde V f) \qquad \mbox{ and } \qquad \tilde V f \,=\, G( \cF, \tilde \cF) (Wf).
\]
By Lemma \ref{semieq} (b) we get
$$
\|f|\TCo Y\| =\|W f |Y\| \asymp \| \tilde  V f |Y \| =\|f| \Co(\tilde \cF,Y)\|.
$$
(Recall that $R$ is assumed to be in $\A_m$ throughout this paper.) 
This implies $\Co(\tilde \cF,Y) = \TCo Y$ with equivalent norms. 
If the canonical dual of $\mathcal{F}$ is $\cA_m$-self-localized then 
it holds  $\Co(\tilde \cF,Y)= \TCo Y = \Co Y$ by Proposition \ref{eqCo}. 
\end{Proof}


Now we study for which class of frames $\mathcal{F}$ the definition of $\Co Y$
does not depend on the particular frame considered.


\begin{proposition}\label{eqFrames}
Assume that $\mathcal{G}=\{g_x\}_{x\in X}$ is an $\cA_m$-self-localized frame 
for $\mathcal{H}$ with  $\cA_m$-self-localized  canonical dual $S_{\cG}^{-1} \mathcal{G}=
\{S_{\cG}^{-1} g_x\}_{x\in X}$, 
where $S_{\cG}$ is the frame operator of $\cG$. If $\mathcal{G} \sim_{\cA_m} \mathcal{F}$ 
for an $\cA_m$-self-localized frame $\mathcal{F}=\{\psi_x\}_{x\in X}$ 
with $\cA_m$-self-localized canonical dual  $S^{-1}_\cF \cF$, then
it holds $\Co(Y)=\Co(\cF,Y) = \Co(\cG,Y)=\Co(S^{-1}_\cF \cF,Y)=\Co(S_{\cG}^{-1} \cG,Y)=\TCo Y$ 
with equivalent norms.
\end{proposition}
\begin{Proof}
By expanding $S^{-1}_\cF \psi_x$ with respect to $S_{\cG}^{-1} \cG$ one has
\[
W_\F f(x) \,=\, \langle f, S^{-1}_\cF \psi_x \rangle \,=\, \int_X \langle f, S_{\cG}^{-1}g_y \rangle 
\langle g_y, S^{-1}_\cF \psi_x  \rangle d\mu(y)
\,=\, G(\cG, S^{-1}_\cF \cF)(W_\cG f)(x),
\]
where $W_\cG f(y) = \langle f, S_{\cG}^{-1} g_y\rangle$.
Lemma \ref{semieq}(b) 
yields $\F \sim_{\cA_m} S^{-1}_\cF \cF$. By Lemma  \ref{semieq}(a)
$G\sim_{\A_m} \F$ and $\F\sim_{\A_m} S^{-1}_\cF \F$ imply
$S^{-1}_\cF \F \sim_{\cA_m} \mathcal{G}$, i.e., $G(\cG, S^{-1}_\cF \cF) \in \cA_m$.
Hence, by $\A_m(Y)\subset Y$
we have $$\|f|\Co(S^{-1} \cF,Y)\| \leq \|G( S^{-1}_\cF \cF,\cG)|\A_m\|\,\|f|\Co(S_{\cG}^{-1} \cG,Y)\|$$ implying
$\Co(S_{\cG}^{-1} \cG,Y) \subset \Co( S^{-1}_\cF \cF,Y)$.
The converse inclusion is shown similarly and with Proposition \ref{eqCo} the proof is completed.
\end{Proof}

\subsection{Intrinsically Localized Duals}\label{sec_dual_loc}

Let $\mathcal{F}=\{\psi_x\}_{x \in X}$ be a continuous frame for $\mathcal{H}$ with (bounded and positive) 
frame operator $S$. One has the following commutative diagram.

\begin{equation}\label{diag}
\begin{matrix}
&\mathcal{H} & \stackrel{S}{\longrightarrow} & \mathcal{H} & \cr &\downarrow V& &\downarrow V& \cr &
\text{ran}(V)&\stackrel{A}{\longrightarrow} &\text{ran}(V)&  \cr &\downarrow \id& \stackrel{A^\dagger}{\nwarrow}&\downarrow \id &
\cr &L^2(X,\mu) & \stackrel{A}{\longrightarrow} & L^2(X,\mu)&
\end{matrix}
\end{equation}

\noindent The operator $A:L^2(X,\mu) \rightarrow L^2(X,\mu)$ can be identified with the
kernel $A(x,y) = G(\cF,\cF)(x,y)= \langle \psi_y,\psi_x\rangle$ and it holds $A = V S W^*$.
By commutativity of the diagram we have $V \circ S = A \circ V$.
Moreover, the operator
$A_{|\text{ran}(V)}: \text{ran}(V) \rightarrow  \text{ran}(V)$
is boundedly invertible and $\text{ker}(A) = \text{ker}(V^*)$.
The operator $A^\dagger = V S^{-1} W^*$ inverts $A$ on $\text{ran}(V)$ and
$\text{ker}(A^\dagger)=\text{ran}(V)^\perp$. Therefore $A^\dagger$ is the (Moore-Penrose) pseudo-inverse of $A$.

\begin{proposition}\label{intrDual}
If $\mathcal{F}$ is a frame  such that $A^\dagger \in \cA$ and $R \in \cA$ then its canonical dual $
S^{-1} \mathcal{F}$ is  intrinsically $\cA$-localized.
\end{proposition}
\begin{Proof}
Since the diagram \eqref{diag} commutes we have
$V S^{-1} = A^\dagger V$.
Applying this equation on $S^{-1}\psi_y$ yields
$G(S^{-1} \cF,S^{-1} \cF)= A^\dagger \circ R$.
With $A^\dagger \in \cA$ and $R \in \cA$ we obtain $S^{-1} \cF \sim_{\cA} S^{-1} \cF$.
\end{Proof}

In the following we will show that any intrinsically $\cA$-localized frame
ensures that $A^\dagger$ and $R$ are in $\cA$, provided $\cA$ is a {\it spectral algebra}.

\begin{definition}\label{spectrAlg} An admissible algebra $\cA$ 
is called a spectral algebra if it fulfills the following additional property
\begin{itemize}
\item[(S)] for all $A=A^* \in \cA$ it holds $\sigma_{\cA}(A)=\sigma(A)$, where $\sigma_{\cA}(A)$ and $\sigma(A)$ are
the spectra of $A$ in $\cA$ and $\cB (L^2(X,\mu) ) $, respectively.
\end{itemize}
\end{definition}
Probably the most important example for our purpose was studied in \cite{BB}.
\begin{examples}\label{A_example}
Assume that $X$ is endowed with a (semi-)metric $d$. Denote
$B_r(x) := \{y \in X, d(x,y) \leq r\}$ the ball of radius $r$ around $x$ and
suppose further that
there exist constants $C,\beta,r_0 \geq 0$ such that $\mu(B_r(x)) \leq C r^\beta$ for all $r \geq r_0$. 
In other words, $X$ is a space of homogeneous type.
Let $\rho:[0,\infty)\to[0,\infty)$ be a concave function
with $\rho(0) = 0$. Then $m(x,y) := e^{\rho(d(x,y))}$ is an admissible weight.
The space $\cA_{2}$ is defined by
\begin{eqnarray*}
\cA_{2}&=&\left\{K : \|K\|_2 =\max\left\{\esssup_{x \in X} \left (\int_{X} |K(x,y)|^2 d\mu(y)\right)^{1/2},
\right.\right.\\
& & \left.\left.
\phantom{\{K : \|K\|_2 =\max\{} \esssup_{y \in X} \left (\int_{X} |K(x,y)|^2 d \mu(x)\right)^{1/2}\right\} < \infty\right\}.\\
\end{eqnarray*}
Endowed with the norm $\| \cdot\|_2$ it is a Banach space. If $\rho$ satisfies the condition
\[
\lim_{\xi \to \infty} \frac{\rho(\xi)}{\xi} \,=\, 0
\]
and
\begin{equation}\label{m_cond}
m(x,y) \,\geq\, (1+d(x,y))^\delta \quad \mbox{for some }\delta >0
\end{equation}
then $\cA_{m,2}= \cA_m \bigcap \cA_{2}$ endowed
with the norm $\|K\|_{m,2}= \max\{\|K|\cA_{m}\|, \|K\|_2\}$ is a spectral algebra. In case of equality in
(\ref{m_cond}) and $0<\delta \leq 1$ this is stated in \cite[Theorem 4.7]{BB}. The general case is proven
completely analogous as in \cite[Theorem 3.1]{gl}.
\end{examples}
For further relevant examples of spectral algebras we refer to \cite{gl}.

The following theorem states that if $A=A^* \in \cA$ for some spectral algebra $\cA$
has a (Moore-Penrose)
pseudo-inverse $A^\dagger$,  then also $A^\dagger \in \cA $. In other words, a
spectral algebra is  ``pseudo-inverse closed''. A proof can be found in  \cite{fg}.
\begin{Theorem}
\label{FoGro}
Let $\cM$ be a closed subspace of $\cH$ with   orthogonal
projection  $P$ onto  $\cM$.   Assume that $A =A^*\in \cA$,
$\text{ker}(A)=\cM^\perp$ and that  $A: \cM \longrightarrow \cM  $ is
invertible.   Then the pseudoinverse $A^\dagger $, i.e., the unique element
in $\cB (\cH )$ satisfying  $A^\dagger A = A A^\dagger = P$ and $
 \text{ker}(A^\dagger)= \cM^\perp$, is an element of $\cA
 $. In particular   $P \in \cA $.
\end{Theorem}


Now let us state the main result of this section.

\begin{Theorem}\label{intrins}
Let $\cA$ be a spectral algebra.
If $\mathcal{F}$ is an  intrinsically $\cA$-localized frame  then also its canonical dual $S^{-1} \mathcal{F}$
is intrinsically $\cA$-localized. In particular, $A^\dagger$ and $R$ are elements of $\cA$.
\end{Theorem}
\begin{Proof}
Since $A=G(\cF,\cF) \in \cA$ is an operator which fulfills the requirement of Theorem \ref{FoGro}, $A^\dagger \in \cA$, $R \in \cA$ and, by Proposition \ref{intrDual},  we obtain $S^{-1}\cF \sim_{\cA} S^{-1}\cF$.
\end{Proof}

\begin{corollary}\label{eqFramesintr} Let $\cA$ be a spectral algebra.
In the set of intrinsically $\cA$-localized continuous frames the relation $\sim_{\cA}$ is an
equivalence and equivalent intrinsically $\cA$-localized continuous frames define 
the same coorbit spaces.
\end{corollary}

We remark that Theorem \ref{intrins} provides a strategy to prove that the kernel $R=G(\F,S^{-1} \F)$
is contained in $\A_m$ (or at least in $\A_1$), which is essential for constructing coorbit spaces and, as
we will see in the following section, also for the extraction of a discrete frame from the continuous one.
In particular examples the appearance of the inverse frame operator $S^{-1}$ in the kernel $R$ makes it hard
to show directly that $R$ is contained in $\A_m$. To overcome this problem 
Theorem \ref{intrins} suggests the following recipe.
One first has to find a suitable subalgebra $\A$ of $\A_m$
which is spectral. (Of course, this is not necessary, if $\A_m$ is spectral itself. 
Unfortunately, it seems
an open question whether $\A_1$ is spectral, see also \cite{gl}.)
For example, an algebra $\A_{m,2}$ as in Example \ref{A_example} might be suitable.
The next step is to prove that $\F$ is $\A$-self-localized, i.e., $G(\F,\F) \in \A$. Potentially, this task 
is much easier since $G(\F,\F)$ does not involve the inverse operator $S^{-1}$. Then 
Theorem \ref{intrins} implies that $R \in \A \subset \A_m$.

\section{Discrete Frames}
\label{sec_5}

In this section we investigate  conditions under which one can extract a discrete 
frame from the continuous one.
In particular, we will derive atomic decompositions and 
{\it Banach frames} for the associated coorbit spaces.

The basic idea is to cover the index set $X$ by some suitable covering $\cU=\{U_i\}_{i \in I}$ with
countable index set $I$ such that the kernel $R$ does not ``vary too much'' on each set $U_i$.
This variation is measured by an auxiliary kernel $\osc_\cU(x,y)$ associated to $R$.
Choosing points $x_i \in U_i$, $i \in I$, we obtain a sampling of the continuous
frame $\{\psi_x\}_{x\in X}$. Under certain conditions on $\osc_\cU$ the sampled system
$\{\psi_{x_i}\}_{i \in I}$ is indeed a frame for $\cH$.
%
%

We start with a definition.
\begin{definition}\label{def_sep} A family ${\cal U} = (U_i)_{i \in I}$ of subsets of $X$
is called 
{\it (discrete) \admissible covering} of $X$ if the following conditions are satisfied.
\begin{itemize}\itemsep-1pt
\item Each set $U_i$, $i \in I$ is relatively compact and has non-void interior.
\item It holds $X = \cup_{i \in I} U_i$.
\item There exists some constant $N>0$ such that
\beq\label{fin_overlap}
\sup_{j\in I} \#\{i \in I, U_i \cap U_j \neq \emptyset\}\leq N < \infty.
\eeq
\end{itemize}
Furthermore, we say that an admissible covering ${\cal U} = (U_i)_{i \in I}$ is {\em moderate} if it fulfills the following additional conditions.
\begin{itemize}
\item There exists some constant $D>0$ such that $\mu(U_i) \geq D$ for all $i\in I$.
\item There exists a constant $\widetilde{C}$ such that
\begin{equation}\label{mu_mod}
\mu(U_i) \,\leq\, \widetilde{C} \mu(U_j) \quad \mbox{for all } i,j \mbox{ with } U_i \cap U_j \neq \emptyset.
\end{equation}
\end{itemize}
\end{definition}

Note that the index set $I$ is countable because $X$ is $\sigma$-compact.
We remark further that we do not require the size of the sets $U_i$ (measured with
$\mu$) to be bounded from above. We only require a lower bound.
Condition (\ref{mu_mod}) means that the sequence $(\mu(U_i))_{i\in I}$ is $\cU$-moderate in the sense of \cite[Definition 3.1]{FG}.
If the sets $U_i$ do not overlap at all, i.e., they form a partition, then this condition is satisfied trivially.
A recipe for the construction of more general \admissible coverings with property \eqref{mu_mod} 
is discussed in \cite{F1} together with some relevant examples. 

For the aim of discretization we have to restrict the class of admissible weight functions (resp. the
class of function spaces $Y$). From now
on we require that there exists a \moderate covering $\cU=(U_i)_{i\in I}$
of $X$ and a constant $C_{m,\cU} $ such that
\beq\label{m_cover}
\sup_{x,y \in U_i} m(x,y) \leq C_{m,\cU} \quad \mbox{ for all } i \in I.
\eeq


Of course, the trivial weight $1$ has this property (provided of course that
\moderate coverings exist),
so that unweighted $L^p(X)$-spaces are admitted. Moreover, if $w$ is a continuous
weight on $X$, then property (\ref{m_cover}) of its associated weight on $X\times X$
defined by (\ref{mw_associate}) means that $w$ is $\cU$-moderate in the terminology
introduced by Feichtinger and Gr\"obner in \cite[Definition 3.1]{FG}.


The next definition will be essential for the discretization problem.

\begin{definition}\label{def_D} A frame $\F$ is said to possess property $D[\delta,m]$
if there
exists a \moderate covering $\cU = {\cal{U}^\delta}=(U_i)_{i \in I}$
of $X$ such that $(\ref{m_cover})$ holds
and such that the kernel $\osc_\cU$ defined by
\[
\osc_\cU(x,y) := \sup_{z \in Q_y} |\langle S^{-1} \psi_x, \psi_y -\psi_z \rangle|
= \sup_{z \in Q_y} |R(x,y) - R(x,z)|,
\]
where $Q_y := \cup_{i, y \in U_i} U_i$, satisfies
\beq\label{osc_delta}
\|\osc_\cU|\A_m\| \,<\, \delta.
\eeq
\end{definition}

We assume from now on that the frame $\F$ possesses at least property $D[\delta,1]$ for some
$\delta>0$. Furthermore,
we only admit weight functions $m$ (resp. spaces $Y$) for which the frame has property
$D[\delta,m]$ for some $\delta>0$.

\subsection{Preparations}

Associated to a function space $Y$ and
to a \moderate covering ${\cal U} = (U_i)_{i\in I}$ we will define two
sequence spaces. Before being able to state their definition we have to make
sure that characteristic functions of compact sets are contained in $Y$.

\begin{lemma}\label{QY} If $Q$ is an arbitrary compact subset of $X$ then the
characteristic function of $Q$ is contained in $Y$.
\end{lemma}
\begin{Proof} Assume that $F$ is a non-zero function in $Y$. Then by solidity we may assume
that $F$ is positive. Clearly, there exists a non-zero
continuous positive kernel $L \in \A_m$. The application of $L$ to $F$ yields a
non-zero positive continuous function
in $Y$ (by the assumption on $\A_m$). Hence, there exists a compact set $U$ with non-void interior
such that $L(F)(x) > 0$ for all $x \in U$. By compactness of $U$ and continuity of $L(F)$
there exists hence a constant $C$
such that $\chi_U(x) \leq C L(F)(x)$ for all $x \in X$. By solidity $\chi_U$ is contained
in $Y$. Now, we set $K(x,y) = \mu(U)^{-1} \chi_Q(x) \chi_U(y)$, which clearly is an element of $\A_m$
by compactness of $Q$ and $U$. It holds $\chi_Q = K(\chi_U)$ and hence $\chi_Q \in Y$.
\end{Proof}

Now we may define the spaces
\begin{align}
Y^\bdiscrete \,:=\, Y^\bdiscrete({\cal U}) \,:=&\, \{(\lambda_i)_{i\in I}, \sum_{i \in I} \lambda_i \chi_{U_i} \in Y\},\notag\\
Y^\ddiscrete \,:=\, Y^\ddiscrete({\cal U})
\,:=&\, \{(\lambda_i)_{i\in I}, \sum_{i \in I} \lambda_i \mu(U_i)^{-1} \chi_{U_i} \in Y\}\notag
\end{align}
with natural norms
\begin{align}
\| (\lambda_i)_{i\in I} | Y^\bdiscrete\| \,:=&\, \|\sum_{i \in I} |\lambda_i| \chi_{U_i}| Y\|,\notag\\
\| (\lambda_i)_{i\in I} | Y^\ddiscrete\| \,:=&\, \|\sum_{i \in I} |\lambda_i| \mu(U_i)^{-1} \chi_{U_i}| Y\|.\notag
\end{align}

If the numbers $\mu(U_i)$ are bounded from above (by assumption they are bounded
from below) then the two
sequence spaces coincide. Lemma~\ref{QY} implies that the finite sequences are contained in 
$Y^\bdiscrete$ and $Y^\ddiscrete$. If the space $(Y,\|\cdot|Y\|)$ is a solid Banach function space, 
then $(Y^\bdiscrete,\|\cdot|Y^\bdiscrete\|)$ and $(Y^\ddiscrete,\|\cdot|Y^\ddiscrete\|)$ are solid BK-spaces, i.e., 
solid Banach spaces of sequences for which convergence implies componentwise convergence 
(this can be seen, for example, as a consequence of Theorem \ref{Yd_prop} (d) and the 
fact that $ Y^\bdiscrete \subset  Y^\ddiscrete$). 
Let us state some further properties of these spaces.

\begin{Theorem}\label{Yd_prop}\begin{itemize}
\item[(a)] The spaces $(Y^\bdiscrete,\|\cdot|Y^\bdiscrete\|)$ and $(Y^\ddiscrete,\|\cdot|Y^\ddiscrete\|)$ are Banach spaces.
\item[(b)] If the bounded functions with compact support are dense in $Y$, then
the finite sequences are dense in $Y^\bdiscrete$ and $Y^\ddiscrete$.
\item[(c)] Denote $a_i := \mu(U_i)$. Further, assume that $w$ is a weight
function on $X$ such that its associated weight $m(x,y) = \max\{w(x)/w(y),w(y)/w(x)\}$ satisfies
(\ref{m_cover}). For $Y = L^p_w(X,\mu),1\leq p \leq \infty$,
it holds $Y^\bdiscrete = \ell^p_{b_p}(I)$ and $Y^\ddiscrete = \ell^p_{d_p}(I)$ with equivalent norms with
\[
b_p(i) \,:=\, a_i^{1/p} \tilde{w}(i), \qquad
d_p(i) \,:=\, a_i^{1/p-1}   \tilde{w}(i)
\]
where $\tilde{w}(i) = \sup_{x \in U_i} w(x)$.
\item[(d)] Suppose that (\ref{m_cover}) holds for the weight function $m$ associated to $Y$ and denote 
$\tilde{v}(i) = \sup_{x\in U_i} v(x)$ and $r(i) = \tilde v(i) \mu(U_i)$. 
Then $Y^\ddiscrete$ is continuously embedded into $\ell^\infty_{1/r}(I)$. 
\end{itemize}
\end{Theorem}

\begin{Proof} The statements (a), (b) and (c) are straightforward to prove. 

For (d) we fix some $k \in I$ and define the kernel
\begin{equation}\label{def_Ki}
K_i(x,y) \,=\, \chi_{U_k}(x) \chi_{U_i}(y),\quad i \in I.
\end{equation} 
For any $i \in I$ we obtain
\[
|\lambda_i| \chi_{U_k} \,=\, K_i( |\lambda_i| \mu(U_i)^{-1} \chi_{U_i})
\,\leq\, K_i(\sum_{j\in I} |\lambda_j| \mu(U_j)^{-1} \chi_{U_j}).
\]
By solidity of $Y$ we get
\begin{align}
|\lambda_i| \|\chi_{U_k}|Y\|
\,&\leq\, \|K_i( \sum_{j\in I} |\lambda_j| \mu(U_j)^{-1} \chi_{U_j})|Y\|
\,\leq\, \|K_i|\A_m\| \, \| \sum_{j\in I} |\lambda_j| \mu(U_j)^{-1} \chi_{U_j} |Y\|\notag\\
&=\, \|K_i|\A_m\| \, \|(\lambda_j)_{j\in I}|Y^\natural\|.\notag
\end{align}
Let us estimate the $\A_m$-norm of $K_i$.
With $y_0 \in U_k$ we obtain
\begin{align}
&\int_X |K_i(x,y)| m(x,y) d\mu(y)
\,\leq\, \chi_{U_k}(x) \int_{U_i} m(x,y) d\mu(y)
\,\leq\, \mu(U_i) \sup_{x \in U_k} \sup_{y \in U_i} m(x,y)\nonumber\\
\,\leq&\, \mu(U_i) \sup_{y\in U_i} m(y,y_0) \sup_{x \in U_k} m(y_0,x)
\,\leq\, C_{m,\cU}  \mu(U_i) \tilde{v}(i).\nonumber
\end{align}
In the last inequality we used that different choices of $z$ in the definition (\ref{v_def})
of $v$ yield equivalent weights. Furthermore, a similar computation yields
\[
\int_X |K_i(x,y)| m(x,y) d\mu(y) \,\leq\, C_{m,\cU}  \mu(U_k) \tilde{v}(i) \leq C_{m,\cU} D^{-1} \mu(U_k) \mu(U_i) \tilde{v(i)}
\]
where $D$ is the constant in Definition \ref{def_sep} of a \moderate covering.
Hence, $\|K_i|\A_m\| \leq C' \mu(U_i) \tilde{v}(i)$ for some suitable constant $C'$ (note that $k$ is fixed). 
This proves the claim.
\end{Proof}

Let us investigate the dependence of the spaces $Y^\bdiscrete$ and $Y^\ddiscrete$
on the particular covering chosen. 

\begin{definition}\label{def_m_equiv} Suppose $\cU = (U_i)_{i \in I}$ and $\cV = (V_i)_{i \in I}$ 
are two \moderate coverings of $X$ over the same index set $I$. Assume that $m$ is a weight function
on $X \times X$. The coverings $\cU$ and $\cV$ are called $m$-equivalent if the following
conditions are satisfied.
\begin{itemize}\itemsep-1pt
\item[(i)] There are constants $C_1,C_2>0$ such that 
$C_1 \mu(U_i) \leq \mu(V_i) \leq C_2 \mu(U_i)$ for all $i \in I$.
\item[(ii)] There exists a constant $C'$ such that $\sup_{x \in U_i} \sup_{y \in V_i} m(x,y) \leq C'$ for
all $i \in I$.
\end{itemize}
\end{definition}

\begin{lemma}\label{lem_mequiv} Let $m$ be the weight function associated to $Y$ and suppose that
$\cU = (U_i)_{i\in I}$ and $\cV = (V_i)_{i \in I}$ are $m$-equivalent \moderate coverings
over the same index set $I$. Then it holds $Y^\bdiscrete(\cU) = Y^\bdiscrete(\cV)$ and 
$Y^\ddiscrete(\cU) = Y^\ddiscrete(\cV)$ with equivalence of norms.
\end{lemma}
\begin{Proof} Assume that $(\lambda_i)_{i\in I}$ is contained in $Y^\bdiscrete(\cV)$. 
Observe that the term 
\[
\int_X \chi_{V_i}(y) \chi_{V_j}(y) d\mu(y) \mu(V_j)^{-1}
\]
equals $1$ for $i=j$ and for fixed $i$ it is non-zero
for at most $N$ different indices $j$ by the finite overlap property (\ref{fin_overlap}).
We obtain
\begin{align}
& \sum_{i \in I} |\lambda_i| \chi_{U_i}(x)
\,\leq\, \sum_{i \in I} |\lambda_i|\sum_{j \in I} \chi_{U_j}(x) 
\int_X \chi_{V_i}(y) \chi_{V_j}(y) d\mu(y) \mu(V_j)^{-1}\notag\\
&=\, \int_X \sum_{i \in I} |\lambda_i| \chi_{V_i}(y) \sum_{j \in I} \chi_{U_j}(x) \chi_{V_j}(y) \mu(V_j)^{-1}
d\mu(y)
\,=\, L(\sum_{i \in I} |\lambda_i| \chi_{V_i})(x),\notag
\end{align}
where the kernel $L$ is defined by
\begin{equation}\label{def_Lkernel}
L(x,y) \,:=\, \sum_{j\in I} \chi_{U_j}(x) \chi_{V_j}(y) \mu(V_j)^{-1}.
\end{equation}
The interchange of summation and integration is always allowed since by the finite overlap
property the sum is always finite for fixed $x,y$.
We claim that $L$ is contained in $\A_m$. 
Using property (i) of $m$-equivalent coverings and once more
the finite overlap property, we get
\begin{align}
\int_X L(x,y) m(x,y) d\mu(y) 
\,&=\, \sum_{j \in I} \chi_{U_j}(x) \int_X \chi_{V_j}(y) \mu(V_j)^{-1} m(x,y) d\mu(y)\notag\\
\,&=\, C' \sum_{j \in I} \chi_{U_j}(x) \,\leq\, C'N \quad \for x \in X.\notag
\end{align}
With property (i) and (ii) in Definition~\ref{def_m_equiv} we get
\begin{align}
&\int_X L(x,y) m(x,y) d\mu(x)
\,=\, \sum_{j \in I} \chi_{V_j}(y) \mu(V_j)^{-1} \int_X \chi_{U_j}(x) m(x,y) d\mu(x)\notag\\
&\leq\, C' \sum_{j \in I} \chi_{V_j}(y) \mu(U_j) \mu(V_j)^{-1}
\,\leq\, C'C_1^{-1}N \quad \for y \in X.\notag
\end{align}
Thus, $L \in \A_m$ and by solidity of $Y$ we conclude that
\[
\|(\lambda_i)_{i\in I})|Y^\bdiscrete(\cU)\| \,\leq\, \|L(\sum_{i\in I} |\lambda_i \chi_{V_i})|Y\|
\,\leq\, \|L|\A_m\|\,\|(\lambda_i)_{i\in}|Y^\bdiscrete(\cV)\|.
\]
Exchanging the roles of $\cU$ and $\cV$ gives a reversed inequality and thus 
$Y^\bdiscrete(\cU) = Y^\bdiscrete(\cV)$. 
Moreover, replacing
$(\lambda_i)_{i\in I}$ by $(\mu(U_i)^{-1} \lambda_i)_{i\in I}$ shows that 
$Y^\ddiscrete(\cU) = Y^\ddiscrete(\cV)$. 
\end{Proof}

For some $i \in I$ we denote $i^*:=\{j \in I, U_i \cap U_j \neq \emptyset\}$.
Clearly, this is a finite set with at most $N$ elements.
The next Lemma states that the sequence spaces $Y^\natural$ are $\cU$-regular in the sense of 
\cite[Definition 2.5]{FG}. 

\begin{lemma}\label{lem_mod} For $(\lambda_i)_{i\in I} \in Y^\natural$ 
let $\lambda_i^+ := \sum_{j \in i^*} \lambda_j$. Then there exists some constant $C>0$ such that 
$\|(\lambda_i^+)_{i\in I}|Y^\natural\| \leq C \|(\lambda_i)_{i\in I}|Y^\natural\|$.
\end{lemma}
\begin{Proof} By Proposition 3.1 in \cite{FG} we have to prove that any permutation $\pi:I \to I$ satisfying
$\pi(i) \subset i^*$ for all $i \in I$ induces a bounded operator on $Y^\natural$, i.e., 
$\|(\lambda_{\pi(i)})_{i\in I}|Y^\natural\| \leq C'\|(\lambda_i)_{i\in I}|Y^\natural\|$.

We define the kernel
\[
K_\pi(x,y) \,:=\, \sum_{i\in I} \mu(U_{\pi^{-1}(i)})^{-1}\chi_{U_{\pi^{-1}(i)}}(x) \chi_{U_i}(y).
\]
It is easy to see that
\[
K_\pi(\mu(U_j)^{-1}\chi_{U_j})(x) \geq \mu(U_{\pi^{-1}(j)})^{-1} \chi_{U_{\pi^{-1}(j)}}(x).
\]
This gives
\begin{align}
\sum_{i\in I} |\lambda_{\pi(i)}| \mu(U_i) \chi_{U_i}(x)\,&=\,
\sum_{i\in I} |\lambda_i| \mu(U_{\pi^{-1}(i)})^{-1} \chi_{U_{\pi^{-1}(i)}}(x)
\,\leq\, \sum_{i\in I} |\lambda_i| K_\pi(\mu(U_i)^{-1} \chi_{U_i})(x)\notag\\
&=\, K_\pi(\sum_{i\in I} |\lambda_i| \mu(U_i)^{-1} \chi_{U_i})(x).\notag
\end{align}
Provided $K_{\pi}$ is contained in $\A_m$ this would give the result by solidity of $Y$.
So let us estimate the $\A_m$-norm of $K_\pi$. We have
\begin{align}
\int_X K_\pi(x,y) m(x,y) d\mu(x) 
\,&=\, \int_X \sum_{i\in I} \mu(U_{\pi^{-1}(i)})^{-1} \chi_{U_{\pi^{-1}(i)}}(x) \chi_{U_i}(y) m(x,y) d\mu(x)\notag\\
&\leq\, \left (\sum_{i\in I} \chi_{U_i}(y) \right )\sup_{i\in I} \sup_{y \in U_i} \sup_{x \in \cup_{j\in i^*} U_j} m(x,y)
\leq C_{m,\cU} ^2 N.\notag
\end{align}
Hereby, we used that for $y \in U_i, x \in U_j$ with $U_i \cap U_j \neq \emptyset$ and $z \in U_i \cap U_j$
it holds $m(x,y) \leq m(x,z) m(z,y) \leq C_{m,\cU} ^2$ by property (\ref{m_cover}). Furthermore by property (\ref{mu_mod}),
we obtain
\begin{align}
&\int_X K_\pi(x,y) m(x,y) d\mu(y) 
\,\leq\, \int_X \sum_{i\in I} \mu(U_{\pi^{-1}(i)})^{-1} \chi_{U_{\pi^{-1}(i)}}(x) \chi_{U_i}(y) m(x,y) d\mu(y)\notag\\
&\leq\, C_{m,\cU} ^2 \sum_{i\in I} \mu(U_{\pi^{-1}(i)})^{-1} \mu(U_i) \chi_{U_{\pi^{-1}(i)}}(y)
\,\leq\, C_{m,\cU} ^2 \widetilde{C} N.\notag
\end{align}
This completes the proof.
\end{Proof}

We will further need a partition of unity (PU) 
associated
to a \moderate covering of $X$, i.e., a family $\Phi = (\phi_{i\in I})_{i\in I}$ of measurable functions that satisfies
$0\leq \phi_i(x) \leq 1$ for all $x\in X$, $\supp \phi_i \subset U_i$
and $\sum_{i\in I} \phi_i(x) = 1$ for all $x \in X$. The construction of such a family $\Phi$ subordinate
to a locally finite covering of some topological space is standard, see also \cite[p.127]{Folland2}.

We may apply a kernel $K$ also to a measure $\nu$ on $X$ by means of
\[
K(\nu)(x) \,=\, \int_{X} K(x,y) d\nu(y).
\]
We define the following space of measures,
\[
D(\cU,M,Y^\ddiscrete) \,:=\, \{\nu \in M_{loc}(X), (|\nu|(U_i))_{i\in I} \in Y^\ddiscrete\}
\]
with norm
\[
\|\nu|D(\cU,M,Y^\ddiscrete)\| \,:=\, \|(|\nu|(U_i))_{i\in I}|Y^\ddiscrete\|,
\]
where $M_{loc}$ denotes the space of 
complex Radon measures.
Spaces of this kind were introduced by Feichtinger
and Gr\"obner in \cite{FG} who called them decomposition spaces.
We identify a function with a measure in the usual way. Then 
\[
D(\cU,L^1,Y^\ddiscrete) \,:=\, \{F \in L^1_{loc}, (\int_{U_i}|F(x)| d\mu(x))_{i\in I} \in Y^\ddiscrete\}
\]
with norm $\|F|D(\cU,L^1,Y^\ddiscrete)\| := \|(\|\chi_{U_i} F|L^1\|)_{i\in I}|Y^\ddiscrete\|$ 
can be considered as a closed
subspace of $D(\cU,M,Y^\ddiscrete)$.

We have the following auxiliary result.

\begin{lemma}\label{lem_Wiener}\begin{itemize}
\item[(a)] It holds $Y \subset D(\cU,L^1,(L^\infty_{1/v})^\ddiscrete)$ with continuous embedding.
\item[(b)] Assume that the frame has property $D[\delta,m]$ for some $\delta > 0$. Then for
$\nu \in D(\cU,M,Y^\ddiscrete)$ it holds $R(\nu) \in Y$ and
$\|R(\nu)|Y\| \leq C \|\nu|D(\cU,M,Y^\ddiscrete)\|$.
\end{itemize}
\end{lemma}
\begin{Proof} (a) Assume $F \in Y$ and let
\[
H(x) := \sum_{i\in I} \|\chi_{U_i} F|L^1\| \mu(U_i)^{-1} \chi_{U_i}(x).
\]
We need to prove $H \in L^\infty_{1/v}$.
Fix $k \in I$. Since $Y$ is continuously embedded into
$L^1_{loc}$ there exists a constant $C$ such that $\|\chi_{U_k} F|L^1\| \leq C \|F|Y\|$
for all $F \in Y$. With $K_i$ as in (\ref{def_Ki}) (and fixed $k\in I$)
it holds $\chi_{U_i} =  \mu(U_k)^{-1} K_i^*(\chi_{U_k})$.
It is shown in the proof of Proposition \ref{Yd_prop} that
$\|K_i|\A_m\| \leq C' \mu(U_i) v(x_i)$ for some constant $C'>0$ and $x_i \in U_i$.
We obtain
\begin{align}
\|\chi_{U_i} F|L^1\| \,=&\, \mu(U_k)^{-1} \|K_i^*(\chi_{U_k}) F|L^1\| \,=\, \mu(U_k)^{-1} \|\chi_{U_k} K_i(F)|L^1\|
\,\leq\, C \mu(U_k)^{-1} \|K_i(F)|Y\| \nonumber\\
\,\leq&\, C \mu(U_k)^{-1} \|K_i|\A_m\|\, \|F|Y\| \leq C'' \mu(U_i) v(x_i) \|F|Y\|. \nonumber
\end{align}
With this we obtain
\[
H(x) \leq C''\|F|Y\| \sum_{i\in I} \chi_{U_i}(x) v(x_i).
\]
For fixed $x$ this is a finite sum over the index set $I_x = \{i \in I, x \in U_i\}$.
It holds
\[
\sup_{i\in I_x} v(x_i) \leq \sup_{i \in I_x} m(x,x_i) m(x,z) \leq C_{m,\cU}  m(x,z) = C_{m,\cU}  v(x)
\]
by (\ref{m_cover}). This proves $H \in L^\infty_{1/v}$ and the embedding is continuous.

(b) Let $\Phi=(\phi_i)_{i\in I}$ be a PU 
associated to ${\cal U}$. Further, we denote
$R_i(x,y):=\phi_i(y) R(x,y)$. Clearly, we have $R(x,y) = \sum_{i \in I} R_i(x,y)$.
We obtain
\[
|R_i(\nu)(x)| \,=\, |\int_X R_i(x,y) d\nu(y)| \,\leq\,
\int_{U_i} |R_i(x,y)| d|\nu|(y) \leq |\nu|(U_i) \|R_i(x,\cdot)\|_\infty.
\]
Observe further that
\[
\mu(U_i) \|R_i(x,\cdot)\|_\infty \, \leq \, \int_X \chi_{U_i}(y) \sup_{z \in U_i} |R(x,z)| dy.
\]
Since the frame is assumed to have property $D[\delta,m]$ we obtain by definition
of $\osc_\cU$ that
\[
|R(x,z)| \leq \osc_\cU(x,y)+|R(x,y)| \quad \mbox{ for all } z,y \in U_i.
\]
This gives
\[
\mu(U_i) \|R_i(x,\cdot)\|_\infty \,\leq \,
\int_X \chi_{U_i}(y) (\osc_\cU(x,y) + |R(x,y)|) dy
\,=\, \left(\osc_\cU + |R|\right)(\chi_{U_i})(x)
\]
and hence,
\begin{align}
\|R(\nu)|Y\|\,=&\, \|\sum_{i\in I} R_i(\nu)|Y\|
\,\leq\, \| \sum_{i\in I} |\nu|(U_i) \mu(U_i)^{-1}(\osc_\cU + |R|)(\chi_{U_i})\|\notag\\
\,=&\,\left\|(\osc_\cU+|R|)\left(\sum_{i \in I}|\nu|(U_i) \mu(U_i)^{-1}\chi_{U_i}\right)|Y\right\|\notag\\
\,\leq&\,(\|\osc_\cU|\A_m\|+\|R|\A_m\|)\,\|\sum_{i \in I}|\nu|(U_i)\mu(U_i)^{-1}\chi_{U_i}|Y\|\notag\\
\label{W_ineq}
\,=&\, (\|\osc_\cU|\A_m\|+\|R|\A_m\|)\,\|\nu|D(\cU,M,Y^\ddiscrete)\|.
\end{align}
This proves the claim.
\end{Proof}

Using this Lemma  
we may prove that the assumption made in Proposition \ref{Prop_proj}
holds in case that the general assumptions of this section are true.

\begin{corollary}\label{cor_RY} If the frame has property $D[\delta,m]$ then $R(Y) \subset L^\infty_{1/v}$ with continuous embedding.
In particular, Proposition~\ref{Prop_proj} holds.
\end{corollary}
\begin{Proof} Suppose $F \in Y$. By Lemma \ref{lem_Wiener}(a) it holds
$F \in D(\cU,L^1,(L^\infty_{1/v})^\ddiscrete)$ and by 
Lemma \ref{lem_Wiener}(b)
we get $R(F) \in L^\infty_{1/v}$.
\end{Proof}

\subsection{Atomic Decompositions and Banach Frames}

Let us give the definition of an atomic decomposition and of a Banach frame. For a Banach
space $B$ we denote its dual by $B^*$.


\begin{definition} A family $\{g_i\}_{i\in I}$ in a Banach space $(B,\|\cdot|B\|)$ is called an atomic decomposition for $B$ if there exist  a BK-space $(B^\ddiscrete(I),\|\cdot|B^\ddiscrete\|)$, $B^\ddiscrete=B^\ddiscrete(I)$, and linear bounded functionals $\{\lambda_i\}_{i\in I} \subset B^*$ (not necessarily unique) such that
\begin{itemize}
\item $(\lambda_i(f))_{i\in I} \in B^\ddiscrete$ for all $f \in B$ and there exists a constant $0<C_1<\infty$ such
that
\[
\|(\lambda_i(f))_{i\in I}|B^\ddiscrete\| \,\leq\, C_1\|f|B\|,
\]
\item if $(\lambda_i)_{i\in I} \in B^\ddiscrete$ then $f = \sum_{i\in I} \lambda_i g_i \in B$ (with unconditional convergence in some suitable topology) and there exists
a constant $0<C_2<\infty$ such that
\[
\|f|B\| \,\leq\, C_2 \|(\lambda_i)_{i\in I}|B^\ddiscrete\|,
\]
\item $f = \sum_{i\in I} \lambda_i(f) g_i$ for all $f \in B$.
\end{itemize}
\end{definition}

We remark that this is not a standard definition (and probably such is not available). 
For instance, Triebel uses this terminology with a slightly different meaning \cite[p.59 and p.160]{Tr2}.
The next definition is due to Gr\"ochenig \cite{gro}.
\begin{definition} Suppose $(B,\|\cdot|B\|)$ is Banach space. A family $\{h_i\}_{i\in I} \subset B^*$
is called a Banach frame for $B$ if there exists a BK-space $(B^\bdiscrete(I),\|\cdot|B^\bdiscrete\|)$, $B^\bdiscrete=B^\bdiscrete(I)$, and a linear bounded reconstruction operator $\Omega:B^\bdiscrete \to B$ such that
\begin{itemize}
\item if $f \in B$ then $(h_i(f))_{i\in I} \in B^\bdiscrete$, and there exist constants $0<C_1,C_2<\infty$
such that
\[
C_1 \|f|B\| \leq \|(h_i(f))_{i \in I}|B^\bdiscrete\| \leq C_2 \|f|B\|,
\]
\item $\Omega(h_i(f))_{i \in I} = f $ for all $f \in B$.
\end{itemize}
\end{definition}

Clearly, these definitions apply also with $B^*$ replaced by the anti-dual $B^\urcorner$.
Now we are prepared to state the main theorem of this article.

\begin{Theorem}\label{thm_atBF} Assume that $m$ is an admissible weight. Suppose the frame $\F=\{\psi_x\}_{x \in X}$ possesses
property $D[\delta,m]$ for some
$\delta$ such that

\begin{equation}\label{cond_delta}
\delta \,(\|R|\A_m\|+\max\{C_{m,\cU}\|R|\A_m\|,\|R|\A_m\|+\delta\}) \,\leq\,1
\end{equation}
where $C_{m,\cU} $ is the constant in (\ref{m_cover}). Let $\cU^\delta$ denote a corresponding \moderate covering of $X$ and choose points $(x_i)_{i\in I} \subset X$ such that
$x_i \in U_i$. Moreover assume that $(Y,\|\cdot|Y\|)$ is a Banach space fulfilling properties (Y1) and (Y2).

Then $\F_d:=\{\psi_{x_i}\}_{i\in I}\subset \K^1_v$ is both an atomic decomposition of $\TCo Y$
with corresponding sequence space $Y^\ddiscrete$ and
a Banach frame for $\Co Y$ with corresponding sequence space $Y^\bdiscrete$.
Moreover, there exists a 'dual frame' $\widehat{\F_d}=\{e_i\}_{i\in I} \subset \H^1_v$
such that
\begin{itemize}
\item[(a)] we have the norm equivalences
\[
\|f|\Co Y\| \,\cong\, \|(\langle f, \psi_{x_i}\rangle)_{i\in I}|Y^\bdiscrete\|
\quad\mbox{ and }\quad \|f|\TCo Y\| \,\cong\, \|(\langle f, e_i\rangle)_{i\in I}|Y^\ddiscrete\|,
\]
\item[(b)] if $f \in \TCo Y$ then
\[
f \,=\, \sum_{i\in I} \langle f,e_i\rangle \psi_{x_i}
\]
with unconditional norm convergence in $\TCo Y$ if the finite sequences are dense in $Y^\ddiscrete$ and with
unconditional convergence in the weak-$*$ topology induced from $(\H^1_v)^\urcorner$ otherwise.
\item[(c)] if the finite sequences are dense in $Y^\bdiscrete$, then for all $f \in \Co Y$ it holds
\[
f \,=\, \sum_{i\in I} \langle f, \psi_{x_i} \rangle e_i
\]
with unconditional convergence in the norm of $\Co Y$.
\end{itemize}
\end{Theorem}

Also discretizations of the canonical dual frame lead to Banach frames and atomic decompositions.

\begin{Theorem}\label{thm_atBF2} Under the same assumptions and with the same notation
as in the previous theorem $\tilde{\F}_d:=\{S^{-1}\psi_{x_i}\}_{ i\in I}\subset \H^1_v$
is both an atomic decomposition of
$\Co Y$ (with corresponding sequence space $Y^\ddiscrete$) and a Banach frame for $\TCo Y$ (with corresponding sequence
space $Y^\bdiscrete$). Moreover, there exists a 'dual frame'
$\widehat{\tilde{\F}_d}= \{\tilde{e_i}\}_{i\in I}\subset \K^1_v$ with the analogous properties
as in the previous theorem.
\end{Theorem}

Let us remark that the two previous theorems hold "uniformly in $Y$". Namely, if $m$ is fixed
then the constant $\delta$ is the same for all function spaces $Y$ satisfying properties (Y1) and (Y2)
with that specific $m$. In particular, the same covering $(U_i)$ can be used for all those spaces $Y$
and $\{\psi_{x_i}\}_{i\in I}$, $x_i \in U_i$, is a Banach frame for all coorbit spaces $\Co Y$ at the same time. 

The previous theorems imply an embedding result.

\begin{corollary} We have the following continuous embeddings
\[
\H^1_v \subset \Co Y \subset (\K^1_v)^\urcorner \qquad\mbox{and}\qquad \K^1_v \subset \TCo Y \subset (\H^1_v)^\urcorner.
\]
\end{corollary}
\begin{Proof} By Proposition \ref{Prop_proj} and Corollary \ref{cor_RY} $f \in \TCo Y$ implies
$Wf \in R(Y) \subset L^\infty_{1/v}$ which in turn means $f \in (\H^1_v)^\urcorner$ by
Lemma \ref{lem_Hdual} and the embedding $\TCo Y \subset (\H^1_v)^\urcorner$
is continuous again by Corollary \ref{cor_RY}. 
Lemma~\ref{QY} shows that
the Dirac element $\delta_i(j) := \delta_{i,j}$ is contained in $Y^\ddiscrete$ and this in turn implies with
Theorem \ref{thm_atBF} that all $\psi_{x_i}$, $i \in I$, are contained in $\TCo Y$ with
$\|\psi_{x_i}|\TCo Y\| \leq \|\delta_i|Y^\ddiscrete\| \leq C \|\delta_i|\ell^1_v\| = C v(x_i)$ .
Since any $x\in X$ may be chosen as one of the $x_i$ it holds $\psi_x \in \TCo Y$ for all $x \in X$
with $\|\psi_x|\TCo Y\| \leq Cv(x)$. Corollary \ref{cor_min} hence implies that $\K^1_v$ is continuously embedded into $\TCo Y$.
The other embeddings are shown analogously.
\end{Proof}

We will split the proof of Theorems \ref{thm_atBF} and \ref{thm_atBF2} into several lemmas. Let us just explain shortly the
idea. Given a \moderate covering $\cU^\delta=(U_i)_{i\in I}$, a corresponding PU $(\phi_i)_{i\in I}$
and points $x_i \in U_i, i\in I$, we define the operator
\[
\UP F(x) \,:=\,\sum_{i\in I} c_i F(x_i) R(x,x_i)\notag
\]
where 
$c_i = \int_X \phi_i(x) d\mu(x)$.
Intuitively, $\UP$ is a discretization of the integral operator
$R$.

If $\UP$ is close enough to the operator $R$ on $R(Y)$ this implies that
$\UP$ is invertible on $R(Y)$ since $R$ is the identity on $R(Y)$ by Proposition \ref{Prop_proj}.
Since $Wf \in R(Y)$ whenever $f \in \TCo Y$ and $R(x,x_i) = W(\psi_{x_i})(x)$ we conclude
\[
W f \,=\, \UP \UP^{-1} W f \,=\, \sum_{i\in I} c_i (\UP^{-1} Wf)(x_i) W \psi_{x_i}
\]
resulting in $f = \sum_{i\in I} c_i (\UP^{-1} Wf)(x_i) \psi_{x_i}$ by the correspondence principle
stated in Proposition \ref{Prop_proj}.
This is an expansion of an arbitrary $f \in \TCo Y$ 
into the elements $\psi_{x_i}, i\in I$, and, thus,
it gives a strong hint that we have in fact an atomic decomposition. Reversing the order of $\UP$ and
$\UP^{-1}$ and replacing $Wf$ by $Vf$
leads to a recovery of an arbitrary $f\in \Co Y$ from its
coefficients $Wf(x_i) = \langle f, \psi_{x_i}\rangle$ and, thus,
we may expect to have a Banach frame. In the following we will make
this rough idea precise. In particular, we need to find conditions on $\delta$ that make
sure that $\UP$ is close enough to the identity on $R(Y)$ (in fact this is ensured by
(\ref{cond_delta})). Moreover, we will need some results that enable us to prove corresponding
norm equivalences.

Let us start with some technical lemmas.

\begin{lemma}\label{lem_syn} Suppose that the frame has property $D[\delta,m]$ for some
$\delta>0$ and that ${\cal U}^\delta=(U_i)_{i\in I}$ is a corresponding \moderate covering of
$X$. Further, assume $(\lambda_i)_{i\in I} \in Y^\ddiscrete$ and $(x_i)_{i\in I}$ to be points
such that $x_i \in U_i$. Then
$x \mapsto \sum_{i\in I} \lambda_i R(x,x_i)$ defines a function in $Y$ and
\beq\label{syn_ineq}
\|\sum_{i\in I} \lambda_i R(\cdot,x_i)|Y\| \,\leq\, C' \|(\lambda_i)_{i\in I}|Y^\ddiscrete\|.
\eeq
The convergence is pointwise, and if the finite sequences are dense in $Y^\ddiscrete$ it is also
in the norm of $Y$.
Furthermore, the series $x \mapsto \sum_{i \in I} R(x,x_i) v(x_i)$ converges pointwise
and absolutely to a function in $L^\infty_{1/v}$.
\end{lemma}

\begin{Proof} Denote by $\epsilon_x$ the Dirac measure in $x$.
Then the application of $R$ to the measure
$\nu:=\sum_{i\in I} \lambda_{i} \epsilon_{x_i}$ results in the function
$x\mapsto \sum_{i\in I} \lambda_i R(x,x_i)$. It follows from Lemma \ref{lem_mod} that
\[
\|\sum_{i\in I} \lambda_{i} \epsilon_{x_i}|D(\cU,M,Y^\ddiscrete)\| \,\leq\, C \|(\lambda_i)_{i\in I}|Y^\ddiscrete\|.
\]
Thus, Lemma \ref{lem_Wiener}(b) yields (\ref{syn_ineq}).
If the finite sequences are dense in $Y$ then clearly the convergence is in the
norm of $Y$.

For the pointwise convergence observe that the space $Y=L^\infty_{1/v}$ has $m$ as
associated weight function. For this choice
it holds $Y^\ddiscrete=\ell^\infty_{1/r}$ where $r(i) = v(x_i)\mu(U_i)$ (Theorem \ref{Yd_prop}(c)).
The application of $|R|$ to
the measure $\nu = \sum_{i\in I} v(x_i) \mu(U_i) \epsilon_{x_i}$ yields $\sum_{i\in I} |R(x,x_i)| v(x_i) \mu(U_i)$.
The estimations 
in (\ref{W_ineq}) are also valid pointwise until the second line, yielding
\[
R(\nu)(x) \leq (\osc_\cU + |R|)(\sum_{i \in I} |\nu|(U_i) \mu(U_i)^{-1} \chi_{U_i})(x)
\]
For our specific choice of $\nu$ we have
$$
|\nu|(U_i) = \sum_{j, x_j \in U_i \cap U_j} |v(x_j)| \mu(U_j) \leq \sum_{j, U_i \cap U_j \neq \emptyset} |v(x_j)| \mu(U_j)< \infty,
$$
since
this is a finite sum. Moreover, for fixed $x$ also 
\[
H(x) = \sum_{i\in I} |\nu|(U_i) \mu(U_i)^{-1} \chi_{U_i}(x).
\]
is a finite sum and, hence, converges pointwise. We already know that $H$ is contained in $L^\infty_{1/v}$.
We conclude that the partial sums of $\sum_{i\in I} |R(x,x_i)| v(x_i) \mu(U_i)$ are dominated
by
\begin{align}
& \int_X (\osc_\cU + |R|)(x,y) H(y) d\mu(y) \,=\, 
\int_X (\osc_\cU + |R|)(x,y) v(y) H(y) v^{-1}(y)d\mu(y) \notag\\
&\leq\, \int_X (\osc_\cU + |R|)(x,y) m(x,y) d\mu(y) m(x,z) 
\sup_{y\in X} |H(y)| v^{-1}(y)\\
&\leq\, m(x,z) (\|\osc_\cU +|R|\,|\A_m\|) \|H| L^\infty_{1/v}\|.\notag
\end{align}
Hence, the sum 
$\sum_{i\in I} |R(x,x_i)| v(x_i) \mu(U_i)$ converges pointwise.
By Theorem \ref{Yd_prop}(d) we have $Y^\ddiscrete \subset \ell^\infty_{1/r}$ for some general $Y$. 
Together with the
results just proven this yields that
the convergence is also pointwise in general.
\end{Proof}

\begin{lemma}\label{lem_BFaux} Suppose that the frame $\F$ has
property $D[\delta,m]$ for some $\delta > 0$
and let ${\cal U}^\delta = (U_i)_{i\in I}$ be an associated \moderate covering of $X$
with corresponding PU $(\phi_i)_{i\in I}$.
If $F \in R(Y)$ then for some constant $D>0$ it holds
\begin{align}\label{ineq_59}
\|\sum_{i\in I} F(x_i) \chi_{U_i}|Y \| \,\leq\, D \| F|Y\|
\quad\mbox{ and }\quad
\|\sum_{i\in I} F(x_i) \phi_i|Y\| \,\leq \, \sigma\|F|Y\|
\end{align}
where $\sigma:=\max\{C_{m,\cU}\|R|\A_m\|,\|R|\A_m\|+\delta\}$
with $C_{m,\cU}$ being the constant in (\ref{m_cover}).
\end{lemma}
\begin{Proof} Since $F \in R(Y)$ it holds $F=R(F)$ by Proposition \ref{Prop_proj} and Corollary \ref{cor_RY}.
This yields
\begin{align}
H(x)\,:=&\, \sum_{i\in I} F(x_i) \chi_{U_i}(x) \,=\, \sum_{i\in I} R(F)(x_i) \chi_{U_i}(x)
\,=\, \sum_{i\in I_x} \int_X R(x_i,y) F(y) \chi_{U_i}(x) d\mu(y)\notag\\
\,=&\, \int_X \sum_{i\in I_x} R(x_i,y) \chi_{U_i}(x) F(y) d\mu(y).\notag
\end{align}
Since the sum is finite over the index set $I_x=\{i, x\in U_i\}$ the interchange of summation
and integration is justified. Defining
\begin{equation}\label{def_K}
K(x,y) \,:=\, \sum_{i\in I} R(x_i,y) \chi_{U_i}(x)
\end{equation}
we obtain $H = K(F)$. We claim that $K \in \A_m$. For the integral with respect to $y$
we obtain
\[
\int_X |K(x,y)| m(x,y) d\mu(y)
\,\leq\, \sum_{i\in I_x} \chi_{U_i}(x) m(x,x_i) \int_X |R(x_i,y)| m(x_i,y) d\mu(y)
\,\leq\, NC_{m,\cU} \|R|\A_m\|
\]
where $N$ is the constant from (\ref{fin_overlap}) and $C_{m,\cU}$ the one from (\ref{m_cover}).
For an estimation of the integral with respect
to $x$ observe first that
\[
|R(x_i,y)| \,\leq\, \osc_\cU^*(x,y) + |R(x,y)|
\]
for all $x\in Q_{x_i} = \cup_{j: U_i \cap U_j \neq \emptyset} U_j$ by definition of $\osc_\cU$.
By Fubini's theorem we obtain
\begin{align}
&\, \int_X |K(x,y)| m(x,y) d\mu(x)
\,=\, \int_{X} \sum_{i\in I} \chi_{U_i}(x) |R(x_i,y)| m(x,y)d\mu(x) \notag\\
\,\leq&\,  \sum_{i \in I} \int_{U_i} (\osc_\cU^*(x,y) + |R(x,y)|) m(x,y) d\mu(x)
\,\leq\, N \int_X (\osc_\cU^*(x,y) + |R(x,y)|) m(x,y) d\mu(x)\notag\\
\,\leq&\, N(\|\osc^*_\cU|\A_m\| + \|R|\A_m\|) \,<\,N(\|R|\A_m\|+\delta)   .\notag
\end{align}
This proves $K \in \A_m$ and we finally obtain
\[
\|\sum_{i\in I} F(x_i) \chi_{U_i}|Y \| \,=\, \|K(F)|Y\| \,\leq\, \|K|\A_m\|\,\|F|Y\|.
\]
A similar analysis shows also the second inequality in (\ref{ineq_59}). 
The constant $N$ from (\ref{fin_overlap}) does not
enter the number $\sigma$ since we replace the characteristic functions by a partition of unity.
\end{Proof}

\begin{corollary}\label{cor_BF} Suppose the frame possesses property $D[\delta,m]$ for some
$\delta > 0$. 
If $f \in \Co Y$ then it holds $\|(V f(x_i))_{i\in I}|Y^\bdiscrete\| \leq C \|f|\Co Y\|$.
\end{corollary}
\begin{Proof}
By Proposition \ref{Prop_proj} it holds $V f \in R(Y)$. By definition of the norm
of $Y^\bdiscrete$ and by Lemma \ref{lem_BFaux} we conclude
$\|(Vf(x_i))_{i\in I}|Y^\bdiscrete\| \leq C \|Vf | Y\| = C \|f|\Co Y\|$.
\end{Proof}

As already announced we need to show that $\UP$ is invertible if $\delta$ is small enough.

\begin{Theorem}\label{thm_Udelta}
Suppose the frame $\F$ possesses property $D[\delta,m]$ for some $\delta >0$.
Then it holds
\begin{equation}\label{Udelta_est}
\|(\id - \UP)|R(Y) \to R(Y)\| \,\leq\, \delta(\|R|\A_m\| +\sigma),
\end{equation}
where $\sigma$ is the constant from Lemma \ref{lem_BFaux}.
Consequently, $\UP$ is bounded and if the right hand side of
(\ref{Udelta_est}) is less or equal to $1$ then $\UP$ is boundedly invertible
on $R(Y)$.
\end{Theorem}
\begin{Proof} Let us first show the implicit assertion that $F \in R(Y)$ implies $\UP(F) \in R(Y)$.
Lemma \ref{lem_BFaux} implies $(F(x_i))_{i\in I} \in Y^\bdiscrete$ which in turns means $(c_i F(x_i))_{i\in I} \in Y^\ddiscrete$.
It follows from Lemma \ref{lem_syn} that $\sum_{i \in I} c_i F(x_i) R(\cdot,x_i)$ converges pointwise to a function
$G=\UP(F) \in Y$. The pointwise convergence implies
the weak-$*$ convergence of $\sum_{i\in I} c_i F(x_i) \psi_{x_i}$ to an element $g$
of $(\H^1_v)^\urcorner$ by Lemma \ref{lem_HK}(b) which is then automatically contained in $\TCo Y$ since $G \in Y$. From
Lemma \ref{lem_HK}(c) follows that $G = Wg = R(Wg)$ and hence $\UP(F) \in R(Y)$.

Let us now introduce the auxiliary operator
\[
\SP F(x) \,:=\, R(\sum_{i\in I} F(x_i) \phi_i)(x).
\]
Assuming $F \in R(Y)$ implies $F=R(F)$ by Proposition \ref{Prop_proj} and Corollary \ref{cor_RY}.
This yields
\[
\|F - \SP F|Y\| \,=\, \|R(F - \sum_{i\in I} F(x_i) \phi_i)|Y\| \,\leq\,
\|R|\A_m\|\,\|F- \sum_{i\in I} F(x_i) \phi_i|Y\|.
\]
We further obtain 
\begin{align}
&\,|F(x)- \sum_{i\in I} F(x_i) \phi_i(x)|
\,=\, |\sum_{i\in I} (R(F)(x)-R(F)(x_i)) \phi_i(x)|\notag\\
\,\leq&\, \sum_{i\in I} \int_X |R(x,y) - R(x_i,y)| |F(y)| d\mu(y) \phi_i(x)
\,\leq\, \sum_{i\in I} \int_X \osc_\cU(y,x) |F(y)| \phi_i(x) d\mu(y)\notag\\
\,=&\, \int_X \osc_\cU^*(x,y) |F(y)| \sum_{i\in I} \phi_i(x) d\mu(y)
\,=\, \osc_\cU^*(F)(x).\notag
\end{align}
Hereby, we used $R(x,y) = \ol{R(y,x)}$, $\supp \phi_i \subset U_i$,
the definition of the kernel $\osc_\cU$
and that $(\phi_i)_{i\in I}$ is a partition of unity. Furthermore, the interchange of
summation and integration in the last line is allowed
since by (\ref{fin_overlap}) the sum is finite
for any fixed $x \in X$. Since $\|K^*|\A_m\| = \|K|\A_m\|$ for all $K \in \A_m$ we
obtain
\begin{equation}\label{Sdelta}
\|F - \SP F|Y\| \,\leq\,\|R|\A_m\|\,\|\osc_\cU^*(F)|Y\|
\,\leq\, \|R|\A_m\|\,\|\osc_\cU|\A_m\|\,\|F|Y\|.
\end{equation}
Let us now estimate the difference of $\UP$ and $\SP$,
\begin{align}
&|\UP F(x) - \SP F(x)| \,=\, |\sum_{i\in I} \int_X \phi_i(y) F(x_i) (R(x,x_i) - R(x,y)) d\mu(y)|\notag\\
\,\leq&\, \sum_{i\in I} \int_X |F(x_i)| \phi_i(y) \osc_\cU(x,y) d\mu(y)
\,=\, \int_X \sum_{i\in I}|F(x_i)| \phi_i(y) \osc_\cU(x,y) d\mu(y).\notag
\end{align}
Denoting $H(y):= \sum_{i\in I}|F(x_i)| \phi_i(y)$ we obtain with Lemma \ref{lem_BFaux} and by solidity of $Y$
\begin{align}
\|\UP F - \SP F|Y\| \,\leq\, \|\osc_\cU(H)|Y\| \,\leq\, \|\osc_\cU|\A_m\|\,\|H|Y\|
\,<\, \delta\sigma \|F|Y\|.
\end{align}
Using the triangle inequality together with \ref{Sdelta} we obtain (\ref{Udelta_est}).
\end{Proof}

Now we have all ingredients to prove Theorem \ref{thm_atBF}.
%
%

\noindent
{\bf\underline{Proof of Theorem \ref{thm_atBF}:}} The condition on $\delta$ implies
by Theorem \ref{thm_Udelta} that $\UP$ is invertible on $R(Y)$.
Assuming $f \in \TCo Y$ means
$Wf \in R(Y)$ by Proposition~\ref{Prop_proj}(a) and Corollary~\ref{cor_RY}. We conclude
\[
W f(x) \,=\, \UP \UP^{-1} Wf(x) =
\sum_{i\in I} c_i\langle \UP^{-1} Wf,\phi_i\rangle R(x,x_i)
=\sum_{i\in I} c_i\langle \UP^{-1} W f,\phi_i\rangle W\psi_{x_i}(x).
\]
Setting $\lambda_i(f) := c_i (\UP^{-1} Wf)(x_i)$ we obtain with
Proposition \ref{Prop_proj}
\begin{equation}\label{f_atomic}
f \,=\, \sum_{i\in I} \lambda_i(f) \psi_{x_i}.
\end{equation}
Since $c_i \leq \mu(U_i)$ we obtain with Lemma \ref{lem_BFaux}
\begin{align}
\|(\lambda_i)_{i\in I}|Y^\ddiscrete\| \,\leq&\, \|(\UP^{-1} Wf)(x_i)|Y^\bdiscrete\|
\,\leq\, C\|\UP^{-1} Wf|R(Y)\|\notag\\
 \,\leq&\, C\on \UP^{-1}|R(Y)\on\,\|f|\TCo Y\|.\notag
\end{align}
Conversely, suppose that $(\lambda_i)_{i\in I} \in Y^\ddiscrete$ and form the function
\[
H(x) \,:=\, \sum_{i\in I} \lambda_i R(x,x_i) \,=\, \sum_{i\in I} \lambda_i W(\psi_{x_i})(x).
\]
Since $Y^\ddiscrete \subset \ell^\infty_{1/\tilde{v}}$ (Theorem \ref{Yd_prop}(d))
the sum converges pointwise to a function in $L^\infty_{1/v}$
by Lemma \ref{lem_syn}. By Lemma~\ref{lem_HK}(b)
the pointwise convergence of the partial
sums of $H$ implies the weak-$*$ convergence in $(\H^1_v)^\urcorner$ of
$f:= \sum_{i\in I} \lambda_i \psi_{x_i}$.
Hence, $f$ is an element of $(\H^1_v)^\urcorner$ and by Lemma~\ref{lem_syn} is therefore contained
in $\TCo Y$. Also from Lemma~\ref{lem_syn} follows
\[
\|f|\TCo Y\| \,=\, \|H|Y\| \,\leq\, C' \|(\lambda_i)_{i\in I}|Y^\ddiscrete\|
\]
and the convergence of the sum representing $f$ is in the norm of $\TCo Y$
if the finite sequences are dense in $Y^\ddiscrete$.
This proves that $\F_d=\{\psi_{x_i}\}_{i\in I}$ is an atomic decomposition of $\TCo Y$.

Now suppose $f \in \Co Y$ and let $F:= V f \in R(Y)$. We obtain
\beq\label{F_inv}
Vf \,=\, \UP^{-1} \UP Vf \,=\, \UP^{-1} \left(\sum_{i\in I} c_i Vf(x_i) W \psi_{x_i}\right).
\eeq
By the correspondence principle (Proposition~\ref{Prop_proj}) this implies
\[
f \,=\, W^* \UP^{-1}  \left(\sum_{i\in I} c_i Vf(x_i) R(\cdot,x_i)\right)
\]
This is a reconstruction of $f$ from the
coefficients $Vf(x_i)=\langle f, \psi_{x_i}\rangle, i\in I$, and the reconstruction
operator $T:Y^\bdiscrete \to \Co Y$,
$T= V^{-1} \UP^{-1} J$ is bounded as the composition of bounded operators.
Note that the operator $J((\lambda_i)_{i\in I})(x) := \sum_{i\in I} c_i \lambda_i R(x,x_i)$
is bounded by
Lemma~\ref{lem_BFaux}.
Setting $Y = L^\infty_{1/v}$ shows that any element of $\Co L^\infty_{1/v} = \K^1_v$ can be
reconstructed in this way.
Now, if for $f \in (\K_v^1)^\urcorner$
it holds 
$(\langle f,\psi_{x_i}\rangle)_{i\in I} \in Y^\bdiscrete$
then the series
$\sum_{i\in I} \langle f, \psi_{x_i} \rangle \phi_i$ converges to an element of $Y$ since
$\phi_i \leq \chi_{U_i}$. 
By bounded invertibility of $\UP$ on $R(Y)$ the right hand side
of (\ref{F_inv}) defines an element in $Y$, hence $f \in \Co Y$.

Using (\ref{F_inv}), the norm equivalence follows from
\begin{align}
&\|f|\Co Y\| \,=\, \|V f|Y\| \,\leq\, \on \UP^{-1}|R(Y) \on \,\|\sum_{i\in I} c_i Vf(x_i)
R(\cdot,x_i)|Y\|\notag\\
\,\leq&\, C\|\UP^{-1}\|\,\|(c_i Vf(x_i))_{i\in I}|Y^\ddiscrete\| \,\leq\, C\on \UP^{-1}\on \,\| (Vf(x_i))_{i\in I}|Y^\bdiscrete\|
\,\leq\, C' \|f|\Co Y\|. \notag
\end{align}
Hereby, we used Lemma~\ref{lem_syn}, $c_i \leq a_i$ and Corollary~\ref{cor_BF}.
Hence, we showed that $\F_d$ is a Banach frame for $\Co Y$.

In order to prove the existence of a dual frame let $E_i := c_i \UP^{-1}(W \psi_{x_i}) \in R(L^1_v)$
and denote $e_i \in \H^1_v$ the unique vector such that $E_i = V(e_i)$. If the finite sequences are
dense in $Y^\bdiscrete$ then
we may conclude from (\ref{F_inv}) by a standard argument (see also \cite[Lemma 5.4]{Gr})
that $f = \sum_{i\in I} \langle f, \psi_{x_i}\rangle e_i$ with norm convergence. This proves (c).

We claim that
\[
\lambda_i(f) \,=\, \langle f,e_i\rangle
\]
yielding together with (\ref{f_atomic}) $f = \sum_{i\in I} \langle f,e_i\rangle \psi_{x_i}$
(with weak-$*$ convergence in general, and if the finite sequences are dense in $Y^\ddiscrete$ with
norm convergence).

If $F \in R(Y)$ then $F(x) = R(F)(x) = \langle F, W \psi_{x} \rangle$. A simple computation shows
\[
\langle \UP F, W \psi_{x} \rangle =
\sum_{i\in I} c_i  F(x_i) \langle R(\cdot,x_i), W\psi_{x}\rangle
=  \sum_{i\in I} c_i F(x_i) \ol{W\psi_{x}(x_i)}
= \langle F, \UP W \psi_{x}\rangle.
\]
Hence, the same relation applies to 
$\UP^{-1}= \sum_{n=0}^\infty (\id - \UP)^n$
and we obtain
\begin{align}
\lambda_i(f) \,=&\, c_i (\UP^{-1} W f)(x_i) \,=\, c_i \langle \UP^{-1} Wf, W \psi_{x_i}\rangle
\,=\, \langle Wf, c_i \UP^{-1} W\psi_{x_i} \rangle\notag\\
\,=&\, \langle Wf, V e_i\rangle \,=\, \langle f, W^* V e_i\rangle \,=\, \langle f, e_i\rangle.\notag
\end{align}
By Lemma \ref{lem_syn} we have the norm estimate
\begin{align}
\|f|\TCo Y\| \,=&\, \|\sum_{i\in I} \langle f,e_i\rangle R(\cdot,x_i)|Y\| \,\leq\, C \|(\langle f,e_i\rangle)_{i\in I}|Y^\ddiscrete\|
\,\leq\, C \|(\UP^{-1}Wf(x_i))_{i\in I}|Y^\bdiscrete\| \notag\\
\,\leq&\,C´\|\UP^{-1} Wf| Y\| \,\leq\, C´\on \UP^{-1}|R(Y) \on \,\|f|\TCo Y\|.\notag
\end{align}
This shows (a) and, thus, we completed the proof of Theorem \ref{thm_atBF}.
Theorem \ref{thm_atBF2} is proved in the same way by exchanging the roles of $V$ and $W$.
{\hspace*{\fill}\qed\vskip1em}

\begin{remark}\label{rem_sample2} 
Using different approximation operators (compare \cite{Gr}) one can prove that under some weaker
condition on $\delta$ one may discretize the continuous frame in order to obtain only atomic decompositions
or only Banach frames with no corresponding results about (discrete) dual frames.
In particular, if 
$\delta\leq1$ then with the procedure
of Theorem~\ref{thm_atBF}
one obtains atomic decompositions and if 
$\delta \leq \|R|\A_m\|^{-1}$
one obtains Banach frames.
\end{remark}

Let us also add some comments about the Hilbert space situation which was the original question of
Ali, Antoine and Gazeau. Here, we need to consider $Y=L^2$ since $\Co L^2 = \TCo L^2 = \H$.
By Lemma~\ref{Yd_prop}(c)
the corresponding sequence space is $Y^\bdiscrete=\ell^2_{\sqrt{a}}(I) = \ell^2(I,a)$ where $a_i = \mu(U_i)$. 
In order to be consistent
with the usual notation of a (discrete) frame it seems suitable to renormalize the frame, i.e.,
under the conditions stated in Theorem~\ref{thm_atBF} (according to Remark~\ref{rem_sample2} 
it is only necessary to have 
$\delta \leq \|R|\A_m\|^{-1}$
 it holds
\[
C_1 \|f | \H\| \leq \sum_{i\in I} |\langle f,\mu(U_i)^{1/2} \psi_{x_i}\rangle|^2 \leq C_2 \|f | \H\|.
\]
This means that $\{\mu(U_i)^{1/2} \psi_{x_i}\}_{ i\in I}$ is a (Hilbert) frame in the usual sense. Of course,
for the aim of Hilbert frames one may choose the trivial weight $m=1$ in Theorem~\ref{thm_atBF}.

One might ask whether the $L^1$-integrability condition $R \in \A_1$ is necessary in order to
obtain a Hilbert frame by discretizing the continuous frame.
The crucial point in the proof of Theorem~\ref{thm_atBF} is that the operator
$\UP$ satisfies
\begin{equation}\label{Udelta_H}
\|\UP - \id| V(\H) \to V(\H)\|<1.
\end{equation}
If one finds a method to prove this
without using integrability assumptions on $R$ then the rest of the proof of Theorem~\ref{thm_atBF} should
still work. However, it is not clear to us how to do this in general.

Concerning a complementary result F\"uhr gave the example of a continuous frame indexed by $\R$ which does
not admit a discretization by any regular grid of $\R$ \cite[Example 1.6.9]{Fuehr}.

\begin{remark} As already noted in Remark \ref{rem_Y2} one may relax condition (Y2) on the function space
$Y$. In this case, one has to restrict to the subalgebra 
$\A:= \A_m \cap \cB(Y)$, where $\cB(Y)$ denotes the continuous
operators on $Y$. The norm on $\A$ is given
by
$
\|K|\A\| \,:=\, \max\{\|K|\A_m\|, \on K|Y\on\}.
$
Note that $\A$ might cease to be closed under the involution $^*$.
In order to carry through all results of this section $Y$ must contain the characteristic functions of
the sets $U_i, i\in I$, which is not automatic if $\A \neq \A_m$. Further, every
occurring kernel must belong to $\A$ and not only to $\A_m$. In particular, $R$ has to be contained in $\A$.
Further, one must replace $\A_m$ by $\A$ in
Definition \ref{def_D} and add that also $\|\osc^*_\cU|Y \to Y\| < \delta$ (which is no longer automatic).
Also in condition (\ref{cond_delta}) on $\delta$ one needs to replace $\A_m$ by $\A$.
Further, one needs to check that
\begin{itemize}
\item $\on K_i|Y\on \leq C v(x_i)$ and $\on K_i^*|Y\on \leq Cv(x_i)$
for the kernel $K_i$ 
defined in (\ref{def_Ki}),
\item the kernel $K(x,y) \,:=\, \sum_{i\in I} R(x_i,y) \chi_{U_i}(x)$ defined in (\ref{def_K}) acts continuously
on $Y$.
\end{itemize}
Under these minor changes one can also invoke the discretization machinery for coorbit spaces
associated to this larger class of function spaces $Y$.
\end{remark}

\section{Localized Continuous Frames Generating Localized Discrete Frames}

In this section we will show that the discretization method presented in the
previous section preserves localization properties.
In particular, we prove that given two localized continuous frames with respect to a suitable
Banach-$*$-algebra $\cA$,  the discretization method 
generates two localized frames with respect to a natural algebra $\A^\bdiscrete$ 
of (infinite) matrices associated to $\cA$,
which is defined similarly as the spaces $Y^\bdiscrete$.

\begin{definition}
Let $\A$ be an admissible Banach-$*$-algebra of kernels on $X \times X$ which is an $\A_m$-bimodule.
Suppose ${\cal U} = (U_i)_{i \in I}$ is a \moderate covering of $X$. Furthermore, for a sequence
$\Lambda = (\lambda_{i,j})_{i,j \in I}$ let
\[
\Lambda_\cU(x,y) \,:=\, \sum_{i,j \in I}  |\lambda_{i,j}| \chi_{U_i}(x) \chi_{U_j}(y).
\]
The algebra $\A^\bdiscrete$ of matrices on $I\times I$ is defined by
\[
\cA^\bdiscrete \,:=\, \cA^\bdiscrete({\cal U})\,= \, \{ \Lambda=(\lambda_{i,j})_{i,j \in I}:
\Lambda_\cU \in \cA\}
\]
with natural norm
\[
\|\Lambda|\A^\bdiscrete\| \,:=\, \|\Lambda_\cU|\A\|.
\]
The multiplication in $\A^\bdiscrete$ is given by
\[
(\Lambda \circ \cE)_{i,j} = \sum_{k \in I} \lambda_{i,k} \epsilon_{k,j} \mu(U_k).
\]
\end{definition}

\begin{proposition}\label{discrA} $\cA^\bdiscrete=\A^\bdiscrete(\cU)$ is a Banach-$*$-algebra.
Moreover, if ${\cal V} = (V_i)_{i \in I}$ is another \moderate
covering of $X$ that is  $m$-equivalent to $\cU$ then $\A^\bdiscrete(\cV) = \A^\bdiscrete(\cU)$
with equivalence of norms. 
\end{proposition}
\begin{Proof}
Let us define $\mu_{k,l} = \int_X \chi_{U_k}(x) \chi_{U_l}(x) d\mu(x)$. Clearly $\mu_{k,k} = \mu(U_k)$.
We have to show that 
$\|\Lambda \circ \cE|\cA^\bdiscrete\|\leq \|\Lambda|\cA^\bdiscrete\| \|\cE|\cA^\bdiscrete\|$,
\begin{align}
&(\Lambda \circ \cE)_\cU(x,y) = \sum_{i,j} \sum_{k \in I} |\lambda_{i,k}| | \epsilon_{k,j}| \mu(U_k)\chi_{U_i}(x)\chi_{U_j}(y)
= \sum_{i,j}\sum_{k \in I} |\lambda_{i,k}| | \epsilon_{k,j}| \mu_{k,k}\chi_{U_i}(x)\chi_{U_j}(y)\notag\\
&\leq \sum_{i,j} \sum_{k,l} |\lambda_{i,k}| | \epsilon_{l,j}| \mu_{k,l} \chi_{U_i}(x)\chi_{U_j}(y)
= \sum_{i,j}\sum_{k,l} |\lambda_{i,k}| | \epsilon_{l,j}| \mu_{k,l}\chi_{U_i}(x)\chi_{U_j}(y)
= (\Lambda_\cU \circ \cE_\cU)(x,y).\notag
\end{align}
We conclude by solidity of $\cA$ and Theorem \ref{Yd_prop}(a) that 
$\cA^\bdiscrete$ is a Banach-$*$-algebra. The second assertion is proven
similarly as in Lemma~\ref{lem_mequiv}. The important
point to note is $\Lambda_\cU \leq L \circ \Lambda_\cV \circ L^*$ with the kernel
$L$ is defined in (\ref{def_Lkernel}).  
\end{Proof}

\begin{proposition}\label{discrAbound}
Let  $\cA$ be an admissible Banach-$*$-algebra  of kernels on $X \times X$ that is an
$\cA_m$-bimodule and satisfies $\cA(Y) \subset Y$. Further, assume
that ${\cal U} = (U_i)_{i \in I}$ is a \moderate covering of 
$X$. Then the algebra $\cA^\bdiscrete({\cal U})$ acts continuously from $Y^\bdiscrete$ into $Y^\bdiscrete$ 
by the mapping
\begin{equation}
\Lambda: Y^\bdiscrete \rightarrow Y^\bdiscrete, \quad
\alpha \mapsto(\Lambda(\alpha))_{i\in I}\,=\,\left(\sum_{j \in I} \lambda_{i,j} \alpha_j \mu(U_j)\right)_{i\in I},
\quad \Lambda = (\lambda_{i,j})_{i,j} \in \A^\bdiscrete.
\end{equation}
In particular, $\A^\bdiscrete$ is continuously embedded into 
$\cB(\ell^2(I,a))$ where $a_i=\mu(U_i)$. 
\end{proposition}
\begin{Proof}
If $\alpha \in Y^\bdiscrete$ then $\alpha_\cU(x) = \sum_{i \in I} \alpha_i \chi_{U_i} \in Y$. A direct computation shows that
$$
\Lambda_\cU(\alpha_\cU)(x) \,=\, \sum_{l,j} \sum_i | \lambda_{l,j}| \alpha_i  \mu_{i,j} \chi_{U_l}(x).
$$
Therefore, by a similar computation as in Proposition \ref{discrA}
\begin{align}
|(\Lambda(\alpha))_\cU(x)| \,\leq\, \sum_{l,j} |\lambda_{l,j}| |\alpha_j|  \mu_{j,j} \chi_{U_l}(x)
\leq \sum_{l,j} \sum_i |\lambda_{l,j}| |y_j| \mu_{i,j} \chi_{U_l}(x)  = \Lambda_\cU(\alpha_\cU)(x).\notag
\end{align}
Since $\cA(Y) \subset Y$ one concludes by solidity of $\cA$ that $\cA^\bdiscrete(Y^\bdiscrete) \subset Y^\bdiscrete$.
In particular, if $Y=L^2(X,\mu)$ then by Theorem \ref{Yd_prop}(c) $Y^\bdiscrete=\ell^2_{a^{1/2}}(I)=\ell^2(I,a)$.
\end{Proof}



Now we have done all preparations to prove the main theorem of this section.

\begin{Theorem}\label{cont2discr}
Let $\cA$ be a Banach-$*$-algebra of kernels on $X \times X$, which is an $\cA_m$-bimodule
with respect to composition. Assume $\cF=\{\psi_{x}\}_{x\in X}$ and $\cG=\{\varphi_x\}_{x\in X}$ 
to be two mutually $\cA$-localized frames
such that $R_{\cF}, R_{\cG} \in \cA_m$. Suppose there exists a \moderate covering
${\cal U^\delta} = (U_i)_{i \in I}$ of $X$ for which both $\cF$ and $\cG$ possess property $D[\delta,m]$ for some $\delta>0$.
Then the discrete systems $\cF_d$ and $\cG_d$, obtained from $\cF$ and $\cG$ 
via $\cF_d=\{\psi_{x_i}\}_{i\in I}$ and $\cG_d=\{\varphi_{y_i}\}_{i\in I}$, 
with $x_i,y_i \in U_i$, are $\cA^\bdiscrete$-localized, i.e., $\cF_d \sim_{\A^\bdiscrete} \cG_d$.
\end{Theorem}
\begin{Proof}
We have to show that $(G(\cF,\cG)(x_i,y_j))_{i,j} \in \cA^\bdiscrete$.
One easily verifies
\[
G(\cF,\cG) \,=\, G(\cF,\cG) \circ R_\cF \qquad \mbox{and} \qquad 
G(\cF,\cG) \,=\, R_\G \circ G(\cF,\cG).
\]
%
Combining these equations 
we obtain
\[
G(\cF,\cG)(x_i,y_j) \,=\,
 \int_X  \int_X  R_{\cG}(x_i,z) R_{\cF}(\xi,y_j)  G(\cF,\cG)(z,\xi) d\mu(\xi)  d\mu(z).
\]
We further deduce
\begin{align}
|G(\cF,\cG)(x_i,y_j)|
\,=&\, \mu(U_i)^{-1} \mu(U_j)^{-1} \int_X \int _X \chi_{U_i}(x) \chi_{U_j}(y)  |G(\cF,\cG)(x_i,y_j)| d\mu(x) d\mu(y)\notag\\
\leq&\, \mu(U_i)^{-1} \mu(U_j)^{-1} \int_X \int _X \chi_{U_i}(x) \chi_{U_j}(y)\int_X  \int_X  \left |R_{\cG}(x_i,z)
R_{\cF}(y_j,\xi)\right|\notag\\
& \phantom{=\, } \times\,|G(\cF,\cG)(z,\xi)| d\mu(\xi)  d\mu(z)d\mu(x) d\mu(y)\notag\\
=&\, \int_X \int _X |G(\cF,\cG)(z,\xi)| \left( \int_{X} |R_{\cF}(y_j,\xi)|  \chi_{U_j}(y) \mu(U_j)^{-1} d\mu(y) \right )\notag\\
& \phantom{=\,} \times \, \left(\int_X  |R_{\cG}(x_i,z)| \chi_{U_i}(x) \mu(U_i)^{-1} d \mu(x)\right )d\mu(\xi)  d\mu(z).\notag
\end{align}
As in the proof of Lemma \ref{lem_Wiener}(b) we have
\begin{align}
& \int_{X} |R_{\cF}(y_j,\xi)|  \chi_{U_j}(y) d\mu(y) \leq (\osc_{\delta}^{\cF} +|R_{\cF}|) (\chi_{U_j})(\xi) \notag\\
\mbox{and } \quad &
\int_X  |R_{\cG}(x_i,z)| \chi_{U_i}(x) d\mu(x) \leq (\osc_{\delta}^{\cG} +|R_{\cG}|) (\chi_{U_i}) (z).\notag
\end{align}
Denoting $\cT_\cF:= \osc_\cU^\cF + |R_\F|\in \A_m$ and $\cT_\cG:=\osc_\cU^\cG + |R_\cG| \in \A_m$ we therefore obtain
\begin{align}
&(G(\cF,\cG)(x_i,y_j))_\cU(x,y)\,=\,\sum_{i,j} |G(\cF,\cG)(x_i,y_j)|\chi_{U_i}(x) \chi_{U_j}(y) \notag\\
\leq&\, \sum_{i,j} \int_X \int _X |G(\cF,\cG)(z,\xi)| \mu(U_i)^{-1} \mu(U_j)^{-1}
\left(\cT_{\cG}(\chi_{U_i})(z) \cT_{\cF}(\chi_{U_j})(\xi) \chi_{U_i}(x) \chi_{U_j}(y) \right)d\mu(\xi)  d\mu(z).\notag
\end{align}
Moreover, with
\[
H^\cG(x,z)\,:=\, \sum_{i\in I}\cT_{\cF}(\chi_{U_i})(z)\chi_{U_i}(x)\mu(U_i)^{-1}, \quad 
H^\cF(y,\xi)\,:=\, \sum_{j\in I}\cT_{\cG}(\chi_{U_j})(\xi) \chi_{U_j}(y)\mu(U_j)^{-1}
\]
we get
\begin{equation}\label{GAb}
((G(\cF,\cG)(x_i,y_j))_{i,j})_\cU \,\leq\, H^\cG \circ |G(\cF,\cG)| \circ (H^\cF)^*.
\end{equation}
Hence, it suffices to show that $H^\cF, H^\cG \in \A_m$. Since for fixed $x$ the sum defining $H^\G(x,z)$
is finite we may interchange the application of $\cT_\G$ with the sum. Denoting
\[
L(x,y) \,:=\, \sum_{i \in I} \chi_{U_i}(x) \chi_{U_i}(y) \mu(U_i)^{-1}
\]
we obtain
\[
(H_\cG)^* \,=\, \cT_\G \circ L \qquad \mbox{and} \qquad (\H^\cF)^* \,=\, \cT_\F \circ L.
\]
Observe that $L$ coincides with the kernel defined in (\ref{def_Lkernel}) for $(V_i)_{i \in I} = (U_i)_{i\in I}$. 
It was already shown in
the proof of Lemma~\ref{lem_mequiv} that $L \in \cA_m$ and hence $H^\cF, H^\cG \in \A_m$.
\end{Proof}

\begin{corollary}\label{locldiscr}
Let $\cA$ be a Banach-$*$-algebra of kernels on $X \times X$ which is an $\cA_m$-bimodule with respect
to composition. Assume $\cF=\{\psi_{x}\}_{x \in X}$  to be an intrinsically $\cA$-localized frame such that there
exists a \moderate covering ${\cal U} = (U_i)_{i \in I}$ of $X$ for which $\cF$ possesses property
$D[\delta,m]$ with $\delta$ small enough (see Theorem~\ref{thm_atBF}). 
Then the discrete system $\cF_d$  generated
from $\cF$  
via $\cF_d=\{\psi_{x_i}\}_{i \in I}$ with $x_i \in U_i$ is an  intrinsically $\cA^{\bdiscrete}$-localized frame for $\cH$.
Moreover, if $\cA^{\bdiscrete}$ is a spectral algebra, then also the canonical
dual of the discrete frame $\tilde \cF_d=\{S_{\cF_d}^{-1} \psi_{x_i}\}_{i \in I}$ is intrinsically $\cA^{\bdiscrete}$-localized, 
where $S_{\cF_d}$ is the frame operator of the discrete frame.
\end{corollary}

In the following example we show that the discretization $A^{\bdiscrete}_{m}$ of an algebra 
$\cA_{m}$ as described in Example \ref{A_example} is again an algebra of this type. In particular,
$\A^\bdiscrete_m$ is then a spectral algebra.

\begin{examples}\label{Ab_example}
Assume that $X$ and $m$ are as in Example~\ref{A_example}. In particular, 
$X$ is endowed with a metric $d$ and $\A_{m,2}$ is a spectral algebra. 
Let $\cU = (U_i)_{i \in I}$ be a \moderate covering satisfying
(\ref{m_cover}). 
Then the discrete algebra $\A_{m}^\bdiscrete=\A_{m,2}^\bdiscrete$ is spectral.
\end{examples}
\begin{Proof} We need to show that $\A^\bdiscrete_m$ is a discrete algebra of the type described
in Example~\ref{A_example}. Let us first note that $A^\bdiscrete_m$ is independent of the choice
of the points $x_i \in U_i$ since by (\ref{m_cover}) 
\begin{equation}\label{m_estim}
m(x_i,x_j) \,\leq\, m(x_i,x) m(x,x_j) \,\leq\, m(x_i,x) m(x_j,y) m(x,y) \,\leq\, C^2 m(x,y)
\end{equation}
for all $x \in U_i, y \in U_j$. Exchanging the roles of $(x,y)$ and $(x_i,x_j)$ gives a reversed inequality.

So let $m^\bdiscrete(i,j):= m(x_i,x_j)$, $i,j \in I$. Clearly, it holds 
$m^\bdiscrete(i,j) = e^{\rho(d^{\,\bdiscrete}\!(i,j))}$ with the (semi-) metric 
$d^{\,\bdiscrete}\!(i,j):=d(x_i,x_j)$ on $I$.
Moreover, $I$ is endowed with the discrete measure $a$ given by $a_i = \mu(U_i)$, $i\in I$.
Denote $B_r^\bdiscrete(i):=\{j: d^{\,\bdiscrete}\!(i,j) \leq r\}$ the ball in $I$ of radius $r$. 
By the finite overlap property (\ref{fin_overlap}) it holds 
\[
a(B_r^\bdiscrete(i)) \,=\, 
\sum_{j \in B_r^\bdiscrete(i)} \mu(U_j) 
\,\leq \, N 
\mu\left( \bigcup_{j \in B_r^\bdiscrete(i)} U_j \right).
\]  
Conditions (\ref{m_cover}) and (\ref{m_cond}) mean that
\[
(1+d(x,y))^\delta \,\leq\, m(x,y) \,\leq\, C_{m,\cU} \quad \for x,y \in U_i,\, i\in I,
\]
for some $\delta > 0$. This implies $d(x,y) \leq C'$ for all $x,y \in U_i$, $i\in I$. 
We conclude that
$
\bigcup_{j \in B_r^\bdiscrete(i)} U_j \subset B_{r+C'}(x_i).
$
Indeed, if $x \in U_j$ with $d(x_i,x_j) \leq r$ then $d(x,x_i) \leq d(x,x_j) + d(x_j,x_i) \leq C'+r$.
Thus, by assumption on the relation of $d$ and $\mu$ (see Example~\ref{A_example}) it holds
$a(B_r(i)) \leq \mu(B_{r+C'}(x_i)) \leq C(r+C')^\beta \leq C'' r^\beta$
for $r \geq r_0'$ with some $r_0'\geq 0$.
Thus, the discrete measure space $(I,a)$ and the weight function $m^\bdiscrete$ satisfy the assumptions
in Example \ref{A_example}. We claim that
\begin{equation}\label{def_Abm_equiv}
\|(\lambda_{i,j})_{i,j}\|'\,:=\,
\max \left\{ \sup_{j \in I} \sum_{i \in I} |\lambda_{i,j}| m^\bdiscrete(i,j) a_i, 
\sup_{i \in I} \sum_{j \in I} |\lambda_{i,j}| m^\bdiscrete(i,j) a_j \right\}
\end{equation}
defines an equivalent norm on $\A^\bdiscrete_m$. 
Indeed for $\Lambda = (\lambda_{i,j})_{i,j} \in \A^\bdiscrete_m$, we obtain using (\ref{m_estim})
\begin{align}
\sup_{j \in I} \sum_{i\in I} |\lambda_{i,j}| m^\bdiscrete(i,j) a_i
\,&=\, \sup_{j \in I} \sum_{i\in I} |\lambda_{i,j}| m^\bdiscrete(x_i,x_j) a_i\notag\\
&=\, \sup_{j\in I} \sum_{i\in I} |\lambda_{i,j}| \int_X \chi_{U_i}(x) m(x_i,x_j) d\mu(x)\notag \\
& \leq \, C_{m,\cU} ^2 
\esssup_{y \in X} \sum_{j\in I} \chi_{U_j}(y) \int_X \sum_{i\in I} |\lambda_{i,j}| \chi_{U_i}(x) m(x,y) d\mu(x)
\notag\\
&=\, C_{m,\cU} ^2 \esssup_{y \in X} \int_X \Lambda_\cU(x,y) m(x,y) d\mu(x).\notag
\end{align}
Using the finite overlap property (\ref{fin_overlap}) and \eqref{m_cover} one similarly obtains the reversed
inequality
\[
 \esssup_{y \in X} \int_X \Lambda_\cU(x,y) m(x,y) d\mu(x)
\,\leq\, NC_{m,\cU} ^2 \sup_{j \in I} \sum_{i \in I} |\lambda_{i,j}| m^\bdiscrete(i,j) a_i.
\]
Exchanging the role of $x$ and $y$ we see that (\ref{def_Abm_equiv}) indeed defines
an equivalent norm on $\A^\bdiscrete_m$. Since $a$ is bounded from below by assumption, i.e, $a_i \geq D$
and $m^\bdiscrete \geq 1$ it holds $\ell^2(I,a) \subset \ell^1_m(I,a)$. This means $\A_2(I,a) \subset \A^\bdiscrete_m$
and $\A^\bdiscrete_m = \A^\bdiscrete_{m,2}$.
\end{Proof}

\begin{remark} In this section we worked with the system $\F_d=\{\psi_{x_i}\}_{i\in I}$. Under 
the assumptions of Theorem~\ref{thm_atBF} this is a discrete frame indexed by $I$ endowed with the measure (weight)
$a$ given by $a_i = \mu(U_i)$. However, for discrete frames one usually prefers to work
with unweighted $\ell^2(I)$-spaces. Indeed, the renormalized frame elements $\F_d^\star=\{\sqrt{a_i} \psi_{x_i}\}_{i\in I}$
form a frame for $\H$ with respect to the unweighted $\ell^2(I)$, see also the end of Section 5. 
Note that the frame operators of 
$\F_d$ and $\F_d^\star$ coincide and, hence, the canonical dual frame of $\F^\star_d$ is a renormalization of the dual
frame of $\F_d$. Also, for the concept of localization it does not play a role
whether one uses $\cF_d$ or $\cF_d^\star$. Indeed, define the map 
$\kappa: \Lambda = (\lambda_{i,j}) \mapsto (\sqrt{a_i a_j}\lambda_{i,j})_{i,j \in I}$ and let 
$\widetilde{\A}^\bdiscrete = \{ \Lambda, \kappa(\Lambda) \in \A^\bdiscrete\}$ with norm 
$\|\Lambda|\widetilde{\A}^\bdiscrete\| = \|\kappa(\Lambda)|\widetilde{\A}^\bdiscrete\|$.
The multiplication in $\widetilde{A}^\bdiscrete$
is defined by
\[
(\Lambda \circ \cE)_{i,j} \,=\, \sum_{k \in I}\lambda_{i,k} \epsilon_{k,j},\quad
\Lambda= (\lambda_{i,j})_{i,j \in I},\, \cE = (\epsilon_{i,j})_{i,j \in I}.
\]
It is easy to see that $\kappa$ is an algebra isomorphism between $\widetilde{\A}^\bdiscrete$ and 
$\A^\bdiscrete$ and, hence, $\widetilde{\A}^\bdiscrete$ is a Banach algebra. $\widetilde{A}^\bdiscrete$ acts
on sequences by $(\Lambda \alpha)_{i \in I} = \sum_{j\in I} \lambda_{i,j} \alpha_j$. 
Moreover, with respect to this action, 
$\widetilde{\A}^\bdiscrete$ is continuously embedded into $\ell^2(I)$ (without weight)
if and only if $\A^\bdiscrete$ is continuously embedded into $\ell^2(I,a)$. Thus,
$\widetilde{\A}^\bdiscrete$ is spectral with respect to $\cB(\ell^2(I))$ if and only if
$\A^\bdiscrete$ is spectral with respect to $\cB(\ell^2(I,a))$.
Now, suppose $\cG_d = \{\phi_{x_i}\}_{i\in I}$ is another discrete frame indexed by $(I,a)$ and denote
by $\cG_d^\star$ its normalization. Then it is easy to see
that $\F_d^\star$ is $\widetilde{\A}^\bdiscrete$-localized with respect to $\cG_d^\star$ if and only if
$\cG_d$ is $\A^\bdiscrete$-localized with respect to $\cF_d$, i.e.,
$ \F_d^\star \sim_{\widetilde{\A}^\bdiscrete} \G_d^\star \Longleftrightarrow \F_d \sim_{\A^\bdiscrete} \G_d$.
\end{remark}

\section{Examples}


\subsection{Classical Coorbit Spaces}
\label{classical}
Of course, the classical theory of Feichtinger and Gr\"ochenig \cite{FG1,FG2,FG3,gro} is a special case of ours.
Let us describe shortly the main features.

Suppose $\G$ is a locally compact, $\sigma$-compact group and $\pi$ an irreducible unitary representation of $\G$
on some Hilbert space. Further assume that $\pi$ is integrable, which means that there exists a non-zero
vector $g \in \H$ such that $\int_{\cG} |\langle g,\pi(x)g\rangle| d\mu(x)<\infty$, where $\mu$ denotes the Haar-measure
of $\G$. This implies that $\pi$ is square-integrable, i.e., there exists a non-zero $g\in \H$ such that
$V_g f \in L^2(\G)$ for all $f \in \H$, where $V_g f(x):= \langle f,\pi(x)g\rangle$ is the (generalized) wavelet transform.
Such a $g$ is called admissible. By a famous theorem of Duflo
and Moore \cite{DM} the space of admissible vectors is dense in $\H$ and it holds
\[
\int_\G |\langle f,\pi(x) g\rangle|^2 d\mu(x) = c_g \|f|\H\|^2.
\]
Thus, $\{\pi(x) g, x\in \G\}$ is a tight continuous frame indexed by $\G$ for any admissible
vector $g$. Since, the frame is tight its frame operator $S$ is a multiple of the identity and, hence, the frame
coincide with its canonical dual (up to normalization). The kernel $R=R_g$ is given by
\[
R_g(x,y) \,=\, \langle \pi(y)g,\pi(x)g\rangle = \langle g, \pi(y^{-1} x)g\rangle = V_g g(y^{-1} x).
\]
Since $\pi$ is assumed to be integrable, $\mu$ is translation invariant and $V_g g(x^{-1}) = \ol{V_g g(x)}$ we
immediately deduce that $R_g$ is contained
in $\A_1$. The application of $R_g$ to a function on $\G$ is a convolution, i.e., $R_g(F) = F*V_g g$.
Thus, it is natural to require the spaces $Y$ to be right $L^1_w$-moduln, i.e., $Y*L^1_w\subset Y$,
where $w$ is a submultiplicative weight function that satisfies some additional assumption, see \cite{FG1,FG2}.
Moreover, one assumes that $Y$ is left and right translation invariant.
If there exists a non-zero $g\in \H$ such that $V_g g$ is contained in $L^1_w$ one may define
the coorbit space $\Co Y$. (Since the frame coincides with its canonical dual we have $\TCo Y = \Co Y$ 
and, thus, it suffices to consider only one class of coorbit spaces.)

For the purpose of discretization one considers discrete admissible coverings of $\G$ of the form $(x_i U)_{i\in I}$
for points $x_i \in \G$ and for some
relatively compact set $U$ with non-void interior. Such coverings exist on every locally compact group.
The condition $\|\osc_\cU|\A_m\|<\delta$ in Definition \ref{def_D} means that $V_g g$ must be contained
in the Wiener amalgam space $W(C_0,L^1_w)$. It is shown in \cite[Lemma 6.1]{FG1} that the set of those
$g$ is dense in $\H$. Furthermore, choosing the set $U$ sufficiently small one can
make $\|\osc_\cU|\A_m\|$ as small as one desires and so with Theorem $\ref{thm_atBF}$ one obtains
atomic decompositions and Banach frames for the corresponding coorbit spaces.

Let us mention some concrete examples.

\noindent
\underline{\bf Homogeneous Besov and Triebel-Lizorkin spaces:} Take $\G= \R^d \rtimes (\R_+^* \times O(d))$, the similitude
group of $\R^d$ with Haar measure $dx a^{-n-1}dadU$. Further, we denote 
$D_a f(t) = a^{-d/2} f(a^{-1}t), a>0$ the dilation, $T_x f(t) = f(t-x),x\in \R^d$ the translation and
$R_U f(t) = f(U^{-1} t), U \in O(d)$ the rotation operator on $L^2(\R^d)$. 
Then $\pi(x,a,U)f := T_x D_a R_U f$ is a square-integrable
irreducible representation of $\G$ on $L^2(\R^d)$ and $V_g f(x,a,U) = \langle f, T_x D_a R_U g\rangle$
is the continuous wavelet transform \cite{SAG}. 

Taking certain mixed norm spaces $L^{p,q}_s$ and tent space $T^{p,q}_s$ \cite{CMS} 
as function spaces $Y$ on $\G$ 
it holds $Y * L^1_{v} \subset Y$.
For suitable Schwartz functions $g$ we have
$V_g g \in L^1_v$ for any of those weight functions $v$, i.e., the kernel $R_g$ is contained
in $\A_m$. Hence, one may define the coorbit $\Co Y$ associated to any of the spaces $L^{p,q}_s$ and
$T^{p,q}_s$.
By a characterization of Triebel in \cite{Tr2} it holds 
\begin{align}
\dot{B}^s_{p,q} = \Co L^{p,q}_{s+d/2-d/q} \qquad \mbox{ and } \qquad \dot{F}^s_{p,q} = \Co T^{p,q}_{s+d/2},\notag
\end{align}
where $\dot{B}_{p,q}^s$ denotes the homogeneous Besov spaces on $\R^d$ and
$\dot{F}_{p,q}^s$ the homogeneous Triebel-Lizorkin spaces.
Theorem \ref{thm_atBF} gives atomic decompositions and Banach frames of wavelet-type for those spaces.
Frazier and Jawerth introduced decompositions of Besov and Triebel-Lizorkin spaces of this type using
the terminology $\phi$-transform \cite{FJ2}. For further details we refer to \cite{gro,hr}.

\underline{\bf Modulation spaces:} 
The original motivation for the construction of the modulation spaces was 
to define Banach spaces of (smooth) functions and distributions 
over a general locally compact 
Abelian (LCA) group $\mathcal{G}$ without having a Lie group structure or a dilation.

%
The (reduced)  Weyl-Heisenberg group $\mathbb{H}_\mathcal{G}$ associated to $\G$ 
is defined as the topological space 
$\mathcal{G}\times \widehat{\mathcal{G}} \times \mathbb{T}$, where $\widehat{\mathcal{G}}$ denotes
the dual group of $\G$ and $\mathbb{T}$ is the torus. The multiplication rule on $\mathbb{H}_{\G}$
is given by 
$$
(x_1,\omega_1,\tau_1) (x_2,\omega_2,\tau_2):=
(x_1+x_2,\omega_1+\omega_2,\tau_1 \tau_2 \omega_1(x_2))
$$
and the Haar measure is the product measure $dx d\omega d \tau$.
The Schr\"odinger representation of $\mathbb{H}_\mathcal{G}$ on 
$\cH=L^2(\mathcal{G})$ is given by 
$\pi(x,\omega,\tau) f(t) := \tau (T_x M_\omega f)(t)$, where $T_x f(t) = f(t-x)$ is the usual translation
and $M_\omega f(t) = \omega(t) f(t)$ is a modulation operator.
Associated to $\pi$ is the short time Fourier transform (STFT) which is  
defined by 
$V_g f (x, \omega, \tau):= \langle f, \pi (x, \omega, \tau)g \rangle = 
\overline{\tau}\langle f, T_x M_\omega g \rangle $ for $f,g \in L^2(\mathcal{G})$.
It is well-known that the Schr\"odinger representation is
indeed square-integrable \cite{Gr} and 
thus $\{\pi (x, \omega, \tau)g\}_{(x, \omega, \tau) \in \mathbb{H}_\mathcal{G}}$ is a 
continuous frame for $L^2(\mathcal{G})$ for any non-zero $g \in L^2(\mathcal{G})$. 

Any coorbit space with respect to the Schr\"odinger representation
of $\mathbb{H}_\mathcal{G}$ is called modulation space.
The most 
prominent examples of modulation spaces are those on the Euclidean space $\R^d$.
Thus, let us assume $\mathcal{G} =\mathbb{R}^d$. We denote $w_s(\omega) = (1+|\omega|)^s$, $s \in \R$, 
a weight function on the frequency variable.
Denote $L^{p,q}_{w_s}(\mathbb{H}_{\mathbb{R}^d})$ the space of measurable 
functions on $\mathbb{H}_{\mathbb{R}^d}$ for which
$$
\|F|L^{p,q}_{w_s}\|:=\left( \int_{\mathbb{R}^d \times \mathbb{T}} \left( \int_{\mathbb{R}^d}   |F(x, \omega,\tau )|^p dx \right )^{q/p}  w_s(\omega)^q d \omega d\tau\right )^{1/q} < \infty.
$$

For $g \in L^2(\mathbb{R}^d)$ such that $V_g g  \in L^1_{w_{|s|}}(\mathbb{H}^d)$, the 
modulation space $M^{p,q}_s(\mathbb{R}^d)$ is defined as the space of tempered 
distributions $f$ such that $V_g f  \in L^{p,q}_{w_s}(\mathbb{H}^{d})$ and hence 
$$
M^{p,q}_{s} \,=\, \Co L^{p,q}_{w_s}.
$$
Furthermore, an application of Theorem \ref{thm_atBF} in this context shows that modulation 
spaces can be characterized by Banach frames of {\it Gabor type}. We refer to \cite{Gr} for further
details and generalizations.

\subsection{\bf Symmetry in Classical Coorbit Spaces} One may also treat
subspaces of the coorbit spaces mentioned above which consist of elements that are invariant under certain
symmetry groups, for instance homogeneous Besov and Triebel-Lizorkin spaces or modulation spaces
of radially symmetric distributions \cite{hr1,hr2}.

Suppose that $\A$ is a compact automorphism group of $\G$ that has also a unitary strongly continuous
representation $\sigma$ on $\H$ such that
\[
\pi(Ax) \sigma(A) = \sigma(A) \pi(x) \quad\mbox{for all }A\in \A, x\in \G.
\]
The space of invariant elements is defined by
\[
\H_\A := \{f \in \H, \sigma(A) f= f \mbox { for all }A\in \A\}.
\]
We denote by $\A x=\{Ax,A\in \A\}$ the orbit of $x$ under $\A$ and define $\K$ to be the space of
all such orbits. $\K$ inherits a natural measure $m$ by projecting the Haar measure $\mu$ of $\G$ onto $\K$.
It is worth to note that $\K$ possesses the structure of a hypergroup.
Further, let
\[
\tilde{\pi}(\A x) = \int_\A \pi(Ax) dA,\quad x\in \G.
\]
The operator $\tilde{\pi}(\A x)$ maps $\H_\A$ into $\H_\A$ for all $\A x \in \K$. Actually, $\pi$ is
a representation of the hypergroup $\K$.
If $\pi$ is square-integrable and $g\in \H_\A\setminus \{0\}$ is admissible
then $\{\tilde{\pi}(\A x)g, \A x \in \K\}$ is a tight continuous frame indexed by $\K$, see \cite{hr1}.
In \cite{hr2} the coorbit theory of these kind of frames is carried
through. In particular, the associated coorbits are subspaces of classical coorbit spaces
consisting of invariant elements. Atomic decompositions and Banach
frames of those spaces could be derived. We remark that here the corresponding sequence spaces
$Y^\bdiscrete$ and $Y^\ddiscrete$ are different from each other in typical situations.
In case of radial modulation spaces, these atomic decompositions
were new. For details we refer to \cite{hr,hr1,hr2}.

We remark that with our results 
one may generalize from the above setting to arbitrary (integrable) representations of
hypergroups. In particular, we expect that the application of our theory to the representations
given in \cite{Roes} leads to a definition of Besov and Triebel-Lizorkin spaces
on Bessel-Kingman hypergroups. These spaces would generalize the radial subspaces of 
Besov-Triebel-Lizorkin spaces
to arbitrary "real-valued dimension". The Besov spaces on the Bessel-Kingman hypergroup
coincide also with the ones
introduced by Betancor {\it et al.} \cite{BMR}.

\subsection{\bf Frames Indexed by Homogeneous Spaces}

As we have discussed in Section \ref{classical}, the group representations
corresponding to the classical integral transforms 
like the wavelet transform and the short time Fourier transform are indeed square-integrable.
However, there are integral transforms related to group representations on $L^2$-spaces on manifolds
which are not square-integrable in a strict sense. In other words, the corresponding 
group is too large.
To overcome this drawback in such cases, an interesting notion of square-integrability 
modulo a subgroup appears in \cite{SAG,SAG1}.

Let $\mathcal{G}$ be a locally compact group and $H$ a closed subgroup of $\G$. Then
the homogeneous space $X=\mathcal{G}/H$ carries a quasi-invariant measure $\mu$.  
Let $\Pi:\mathcal{G} \rightarrow X$ denote the canonical projection. Moreover, suppose
that $\sigma:X \rightarrow \mathcal{G}$ is a measurable section of $\mathcal{G}$, i.e., $\Pi \circ \sigma (x) = x$, 
for all $x \in X$. We say that a unitary representation $\pi$ of $\mathcal{G}$ on $\mathcal{H}$ 
is square-integrable $\text{mod}(H,\sigma)$ provided there exists some $g \in  \mathcal{H}$ 
such that $\{\pi(\sigma(x)) g\}_{x \in X}$ is a continuous frame for $\mathcal{H}$.
Many important examples can be described in this setting,
such as the continuous wavelet transform on the sphere introduced by Antoine and Vandergheynst \cite{AV} 
and a notion of Gabor transform on the sphere developed by Torresani \cite{T3}.
Also a mixture of Gabor and wavelet transform on $\R^d$ fits into this approach \cite{HL2,T1} (see also
the example on $\alpha$-modulation spaces).

As a matter of fact, the theory of Feichtinger and Gr\"ochenig is no longer applicable in this setting.
Efforts to adapt the original coorbit space theory to homogeneous spaces have 
been done recently by Dahlke {\it et al.} \cite{DST,DST1}, allowing for instance the definition of 
{modulation spaces} on spheres as coorbit spaces. However,  their approach 
works under the assumption that the frame is tight. This fails to be true for the continuous wavelet transform on the
sphere \cite{AV} and the mixed Gabor / wavelet transform \cite{HL2,T1}. 
Clearly, in the present paper we avoided this drawback by permitting general continuous frames. 
However in order to apply our results to the mentioned examples, still some effort has to be done.
In particular, one needs to check that the corresponding kernel $R$ is contained in $\A_m$ and
that there exists a suitable covering $\cU$ of the index set $\G/H$ such that the resulting kernel $\osc_\cU$
satisfies $\|\osc_\cU|\A_m\| < \delta$ for some $\delta$, which is small enough.
It seems that this task is rather difficult for the examples mentioned above. We postpone
detailed discussions to later contributions.

\subsection{Non-standard Examples}

In this subsection we collect two relevant examples where neither classical coorbit space theory
\cite{FG1,FG2,FG3,gro} nor its  recent generalizations \cite{DST,DST1,hr2}
can be applied.

\noindent
\underline{\bf Inhomogeneous Besov and Triebel-Lizorkin spaces:}

Suppose that
$\psi$ is a radial Schwartz
function on $\R^d$ with $\supp \hat{\psi} \subset \{x, 1/2 \leq |x| \leq 2\}$ such that
\[
\int_{\R^d} \frac{|\hat{\psi}(x)|^2}{|x|^d} dx \,=\, c_\psi = 1,
\]
i.e., $\psi$ is an admissible wavelet for the continuous wavelet transform on $\R^d$. Hereby, $\hat{\psi}$
denotes the Fourier transform of $\psi$ and $|x|$ the Euclidean norm of $x \in \R^d$.
Further, let $\phi$ be a Schwartz function on $\R^d$ such that
\begin{equation}\label{wave_cond2}
|\hat{\phi}(\xi)|^2 + \int_0^1 |\hat{\psi}(t\xi)|^2 \frac{dt}{t} = 1 \quad\mbox{ for all }\xi \in \R^d.
\end{equation}
Observe that $1=c_\psi = \int_0^\infty |\hat{\psi}(t\xi)|^2 dt/t$ and hence $0\leq |\hat{\phi}(\xi)|^2 \leq 1$.
The support condition on $\hat{\psi}$ implies that $|\hat{\phi}(x)|^2 = 1$ for all $|x|\leq 1/2$
and $\hat{\phi}(x) = 0$ for all $|x|\geq 2$. With the unitary dilation $D_t$ and translation $T_x$
define
\begin{align}
\psi_{\infty, x}(y) \,=&\, T_x \phi(y) \,=\, \phi(y-x),\notag\\
\psi_{t,x}(y)       \,=&\, T_x D_t \psi(y) \,=\, t^{-d/2}\psi(t^{-1}(y-x)).\notag
\end{align}
(Here, $\infty$ denotes a separated point.)
An straightforward computation shows that as a consequence of (\ref{wave_cond2}) it holds
\begin{equation}\label{frame_eqpsi}
\int_{\R^d} \left(|\langle f, \psi_{\infty,x} \rangle|^2 + \int_0^1 |\langle f, \psi_{t,x}\rangle|^2
\frac{dt}{t^{d+1}}\right) dx \,=\, \|f\|^2 \quad \mbox{for all } f \in L^2(\R^d).
\end{equation}
We set $X:= (\{\infty\} \cup (0,1]) \times \R^d$. Then (\ref{frame_eqpsi}) means that
$\{\psi_{t,x}, (t,x) \in X\}$ is a continuous frame indexed by $X$ with associated measure
\[
\int_X F(z) d\mu(z) \,=\, \int_{\R^d} F(\infty,x) + \int_0^1 F(t,x) \frac{dt}{t^{d+1}} dx.
\]
Since this frame is tight the associated frame operator is the identity and the frame coincides with its canonical
dual. We remark that the index set $X$ apparently does not have the structure of a group, of a homogeneous space
or of an orbit space.
Further, for $s\in \R$ we define $w_s(t,x) := t^s$ for $t \in (0,1]$ and $w_s(\infty,x) := 1$.
Its associated weight $m_s$ (\ref{mw_associate}) becomes
\[
m_s((t,x),(r,y)) \,=\, m_s(t,r) \,=\, \left(\max \left\{\frac{t}{r},\frac{r}{t}\right\}\right)^{|s|}
\]
for $r,t \neq \infty$ (with obvious modification if $r=\infty$ or $t=\infty$).
It is straightforward to show that the kernel
\[
R((t,x),(r,y)) \,=\, \langle \psi_{r,y},\psi_{t,x} \rangle
\]
is contained in $\A_{m_s}$ for all $s\in \R$. With $L^p_s := L^p_{w_s}$ 
it holds $\A_{m_s}(L^p_s) \subset L^p_s$
and, thus, we may define the coorbit spaces $\Co L^p_{m_s}$, $1\leq p \leq \infty$. 
A characterization of Triebel 
in \cite[Sections 2.4.5 and 2.5.3]{Tr2} shows
\begin{equation}\label{Co_inhomog}
B_{p,p}^s(\R^d) \,=\, \Co L^{p}_{s+d/2-d/p},
\end{equation}
where $B_{p,p}^s$ is an inhomogeneous Besov space. Note that $B_{p,p}^s=F_{p,p}^s$, where the latter denotes
an inhomogeneous Triebel-Lizorkin space. 
We remark that one can extend (\ref{Co_inhomog}) also to Besov and Triebel-Lizorkin spaces with $p \neq q$.
To do this one needs to introduce mixed norm spaces $L^{p,q}_s$ and tent spaces $T^{p,q}_s$.
According to Remark~\ref{rem_Y2} one has to check that $R(L^{p,q}_s) \subset L^{p,q}_s$ 
and $R(T^{p,q}_s) \subset T^{p,q}_s$. Note that this does not follow automatically
from $R \in \A_m$. 

Our discretization machinery (Theorem~\ref{thm_atBF}) yields wavelet type Banach 
frames and atomic decompositions
of the inhomogeneous Besov and Triebel-Lizorkin spaces, similarly as in \cite{FJ2}.
We postpone the details eventually to successive contributions.

\underline{\bf $\alpha$-modulation spaces:} 
Consider a system of functions of the type
\begin{equation}
\label{aframe}
 \mathcal{G}_\alpha(g):=\{M_{\omega} D_{(1+|\omega|)^{-\alpha}} T_x g\}_{x,\omega \in \mathbb{R}},
\end{equation}
where $\alpha \in [0,1)$ and $g \in L^2(\mathbb{R})$. 
One has the following 
well-known result (see for example \cite{HN}).
\begin{proposition}
\label{acontframe}
Let us fix $\alpha \in [0,1)$. Suppose $g \in L^2(\mathbb{R})$ such that there exists a constant $A>0$ 
for which the function
$$
        \sigma_g^\alpha(\xi) \,:=\,  \int_{\mathbb{R}} 
\left |\widehat{g} \left (\frac{\xi -w}{(1+|w|)^\alpha}\right) \right|^2  (1+|w|)^{-\alpha} dw
$$
satisfies
\begin{equation}\label{a_admiss}
        A^{-1} \leq \sigma_g^\alpha(\xi) \leq A, \quad \text{ for almost all } \xi \in \mathbb{R}.
\end{equation}
Then $\mathcal{G}_\alpha(g)$ is a continuous frame for $L^2(\mathbb{R})$. 
A typical function satisfying \eqref{a_admiss} for all $\alpha \in [0,1)$ is the Gaussian.
\end{proposition}
Associated to such frames is the {\it $\alpha$-transform}
\begin{equation}
\label{atrans}
                V_g^\alpha f (x,\omega)= \langle f, M_\omega  D_{(1+|\omega|)^{-\alpha}} T_x g \rangle, \quad (x,\omega) \in \mathbb{R}^2.
\end{equation} 
One easily verifies
that for $\alpha=0$ the family \eqref{aframe} is in fact a continuous Gabor frame and $V_g^0$ is the short time Fourier transform, while for $\alpha \rightarrow 1$ the family tends to the situation encountered in the wavelet context, i.e.,
$V_g^1$ is a slight modification of the classical wavelet transform.
The intermediate case $\alpha=1/2$ gives the Fourier-Bros-Iagolnitzer transform \cite{BI}.

In \cite{HN,HL2,ff,f1} a new class of spaces has been suggested as retract of weighted $L^{p,q}$ spaces by means 
of $V_g^\alpha$, in the same way as it is done for the more classical cases of modulation and Besov spaces. 
In particular, in \cite[Theorem 3.5]{f1} it has been shown that this class 
coincides with the family of so called $\alpha$-modulation spaces $M^{p,q}_{s,\alpha}$ introduced independently 
by Gr\"obner \cite{G} and 
P\"aiv\"arinta/Somersalo \cite{PS} as an ``intermediate'' family between modulation and {\it inhomogeneous} Besov spaces.
In particular, it holds 
\begin{equation}
M^{p,q}_{s+\alpha(1/q-1/2),\alpha}(\mathbb{R})=\{f \in \mathcal{S}'(\mathbb{R}): V_g^\alpha(f) \in L^{p,q}_{w_s} (\mathbb{R}^2)\},
\end{equation} 
\begin{equation} \|f |  M^{p,q}_{s+\alpha(1/q-1/2),\alpha}\| \asymp  \left(\int_{\mathbb{R}} \left( \int_{\mathbb{R}} |V_g^\alpha(f)(x,\omega)|^p dx \right)^{q/p} (1+|\omega|)^{s q} d\omega \right )^{1/q}
\end{equation}
where $g$ is a suitable Schwartz function and $L^{p,q}_{w_s}(\mathbb{R}^2)$ is the space of 
functions $F$ on $\mathbb{R}^2$ such that
$$
        \|F|L^{p,q}_{w_s}\|:= \left(\int_{\mathbb{R}} \left( \int_{\mathbb{R}} |F(x,\omega)|^p dx \right)^{q/p} (1+|\omega|)^{s q} d\omega \right )^{1/q} < \infty.
$$
For $\alpha=0$ the space $M^{p,q}_{s,0}(\mathbb{R})$ coincides  with the {modulation space} $M^{p,q}_s(\mathbb{R})$ 
and  for $\alpha \rightarrow 1$  we obtain 
the {inhomogeneous Besov space} 
$B^{p,q}_{s}(\mathbb{R})=M^{p,q}_{s,1}(\mathbb{R}) $. 
The $\alpha$-modulation spaces are known to have nice analysis properties. 
For instance, the mapping properties of 
pseudodifferential operators in H\"ormander classes on $\alpha$-modulation spaces as by  
Holschneider, Nazaret \cite{HN} and  Borup \cite{B} generalize classical results of Cordoba and 
Fefferman \cite{CF}. Moreover, we expect that such spaces have a key role in the study of pseudodifferential 
operators modeling the transmission of (digital) signals in wireless communications and in corresponding 
numerical methods  \cite{DFR}.

The description of $\alpha$-modulation spaces as coorbit spaces associated to the continuous frame 
$\mathcal{G}_\alpha(g)$ is still a matter of investigation. In fact, while the square-integrability 
of $ V_g^\alpha(f)(x,\omega)$ is ensured by Proposition \ref{acontframe}, the localization properties of 
the corresponding kernel 
$K(x,\omega;\tilde x, \tilde \omega):=\langle M_{\omega} D_{(1+|\omega|)^{-\alpha}} T_x g, M_{\tilde \omega} D_{(1+|\tilde \omega|)^{-\alpha}} T_{\tilde x} g \rangle$ have not yet been proven for $\alpha \in (0,1)$. 

\subsection{Reproducing Kernel Hilbert Spaces and Sampling}

Let $\H \subset L^2(X,\mu)$ be a reproducing kernel Hilbert space with reproducing kernel $K_x(t)$, $t,x \in X$.
In particular, we have $f(x) = \langle f, K_x\rangle$ for all $f \in \H$. This gives
\[
\int_X |\langle f, K_x\rangle|^2 d\mu(x) \,=\, \int_X |f(x)|^2 d\mu(x) \,=\, \|f|\H\|^2 \quad \mbox{for all } f\in \H. 
\]
Hence, the family $\{K_x\}_{x \in X}$ is a tight continuous frame for $L^2(X,\mu)$ 
with frame operator $S$ being
the identity.
The corresponding kernel $R$ is given by
\[
R(x,y) \,=\, \langle K_y,K_x\rangle \,=\, K_y(x),\quad x,y \in X.
\]
Moreover, the transform $V$ is the identity on $\H$, i.e., $V f(x) \,=\, \langle f, K_x\rangle = f(x)$.
Let $w$ be some weight function on $X$ and $m$ its associated weight (\ref{mw_associate}).
Provided $R$ is contained in $\A_m$ then the coorbit spaces are well-defined and we have
\[
\Co L^p_w \,=\, \left\{f \in L^p_w, f(x) = \int_X f(y) \ol{K_x(y)} d\mu(y)\right\}.
\]
If the continuous frame $\{K_x\}_{x \in X}$ possesses 
property $D[\delta,m]$ for some $\delta$ small enough, then we may envoke the discretization
machinery. This yields sampling theorems for $\Co L^p_w$, in particular for $\H = \Co L^2$. 
Indeed, if $\{K_{x_i}\}_{i \in I}$ forms a Banach frame then
it holds
\[
f(t) \,=\, \sum_{i\in I} \langle f, K_{x_i} \rangle e_i(t)\,=\, \sum_{i\in I} f(x_i) e_i(t),\quad t\in X,
\]
for all $f \in \Co L^p_w$. Hereby, the functions $\{e_i\}_{i\in I}$ form a dual frame in the sense of 
Theorem~\ref{thm_atBF}. Moreover, if $\{K_{x_i}\}_{i\in I}$ constitutes an 
atomic decomposition then we have an expansion
\[
 f(t) \,=\, \sum_{i\in I} \langle f, e_i\rangle K_{x_i}(t), \quad t \in X,
\]
for all $f \in \Co L^p_w$, $1 \leq p < \infty$.
{
\noindent
{\bf Massimo Fornasier}, Dipartimento di Metodi e Modelli Matematici per le Scienze Applicate, Universit\`a di Roma ``La Sapienza'', Via A. Scarpa, 16/B, I-00161 Roma, Italy;\\
NuHAG (Numerical Harmonic Analysis Group), Fakult\"at f\"ur Mathematik, Universit\"at Wien, Nordbergstrasse, 15, A-1090 Wien, Austria.\\ email: \underline{mfornasi@math.unipd.it}\\

\noindent {\bf Holger Rauhut}, Zentrum Mathematik, Technische Universit\"at M\"unchen,  Bolzmannstr. 3, D-85747 Garching, Germany. \\email: \underline{rauhut@ma.tum.de}

}

\begin{thebibliography}{9}
\bibitem{SAG} {S. T. Ali, J. P. Antoine, J. P. Gazeau}, \textit{Coherent States, Wavelets and their Generalizations}, Springer-Verlag, 2000.
\bibitem{SAG1} {S. T.  Ali, J. P. Antoine, J. P. Gazeau}, {Continuous frames in Hilbert spaces}, Annals of Physics, {222}, 1--37, 1993.
\bibitem{AV}  J. P. Antoine, P. Vandergheynst, Wavelets on the 2-sphere: A group-theoretical approach, Appl. Comp. Harm.
Anal., 7, 262--291, 1999.
\bibitem{BB} B. A. Barnes, The spectrum of integral operators on Lebesgue spaces, J. Oper. Theory, 18, 115--132, 1987.
\bibitem{BMR} J.J. Betancor, L. Rodr\'iguez-Mesa, {On the Besov-Hankel spaces}, J. Math. Soc. Japan, 50, 781--788, 1998.
\bibitem{Bonsall} F.F. Bonsall, {Decomposition of functions as sums of elementary functions}, Quart. J. Math. Oxford, 37(2), 
129--136, 1986.
\bibitem{B} L. Borup, Pseudodifferential operators on $\alpha$-modulation spaces,  to appear in J. Func. Spaces and Appl., {2}(2), 2004.
\bibitem{BI} J. Bros, D. Iagolnitzer, Support essentiel et structure analytique des distributions, in Seminaire Goulaouic-Lions-Schwartz, exp. no. 18, 1975.
\bibitem{CMS} R.R. Coifman, Y. Meyer, E.M. Stein, {Some new function spaces and their applications to harmonic analysis},
J. Funct. Anal., {62}, 304--335, 1985.
\bibitem{CF} A. Cordoba, and C. Fefferman, Wave packets and Fourier integral operators, Comm. Part. Diff. Eq., 3, 979--1005, 1978.
\bibitem{DFR} S. Dahlke, M. Fornasier, and T. Raasch,  Adaptive frame methods for elliptic operator equations, Bericht Nr. 2004-3, Fachbereich Mathematik und Informatik, Philipps-Universit\"at Marburg, 2004.
\bibitem{DST} S. Dahlke, G. Steidl, and G. Teschke, Coorbit spaces and Banach frames on homogeneous spaces with applications to the sphere, Adv. Comp. Math., 21, 147--180, 2004. 
\bibitem{DST1} S. Dahlke, G. Steidl, and G. Teschke, Weighted coorbit spaces and Banach frames on homogeneous spaces, preprint, Bericht Nr. 2003-4,  Fachbereich Mathematik und Informatik, Philipps-Universit\"at Marburg, 2003, to appear in J. Four. Anal. Appl.
\bibitem{D2}  {I. Daubechies}, {\it Ten Lectures on Wavelets}, SIAM, 1992.
\bibitem{DG} I. Daubechies, A. Grossmann, Y. Meyer, Painless nonorthogonal expansions, J. Math Phys., 27(5), 1271--1283, 1986.
\bibitem{DS}  {R. J. Duffin, A. C. Schaeffer}, {A class of nonharmonic Fourier series}, Trans. Amer. Math. Soc., {72}, 341--366, 1952.
\bibitem{DM} M. Duflo, C.C. Moore, {On the regular representation of a nonunimodular locally compact group},
J. Funct. Anal., {21}, 209--243, 1976.

\bibitem{F1} {H. G. Feichtinger}, {Banach spaces of distributions defined by decomposition methods II}, Math. Nachr., 132, 207--237, 1987.
\bibitem{F4} {H. G. Feichtinger}, {Atomic characterization of modulation spaces through Gabor-type representations}, Proc. Conf. Constr. Function Theory, Rocky Mountain J. Math., 19, 113--126, 1989.
\bibitem{ff} H. G. Feichtinger, M. Fornasier, Flexible Gabor-wavelet atomic decompositions for $L^2$-Sobolev spaces, to appear in Annali di Matematica Pura e Applicata.
\bibitem{FG} {H. G. Feichtinger, P. Gr\"obner}, {Banach spaces of distributions defined by decomposition methods I},  Math. Nachr., {123}, 97--120, 1985.
\bibitem{FG1} {H. G. Feichtinger, K. Gr\"ochenig}, {A unified approach to atomic decomposition via integrable group representations}, in ``Proc. Conference on Functions, Spaces and Applications, Lund, 1986'', Springer Lect. Notes Math., {1302}, 52--73, 1988.
\bibitem{FG2}  {H. G. Feichtinger, K. Gr\"ochenig}, {Banach spaces related to integrable group representations and their atomic decompositions I}, J. Funct. Anal., {86}, 307--340, 1989.
\bibitem{FG3} {H. G. Feichtinger, K. Gr\"ochenig}, {Banach spaces related to integrable group representations and their atomic decompositions II}, Monatsh. Math., {108}, 129--148, 1989.
\bibitem{Folland2} G. B. Folland, {\it Real Analysis}, John Wiley \& Sons, New York, 1984.
\bibitem{f1} M. Fornasier, Banach frames for $\alpha$-modulation spaces, Preprint 10/2004, Dipartimento di Metodi e Modelli Matematici per le Scienze Applicate, Universit\`a di Roma ``La Sapienza'', 2004.
\bibitem{fg} M. Fornasier, K. Gr\"ochenig, Intrinsic localization of frames, Preprint 8/2004, Dipartimento di Metodi e Modelli Matematici per le Scienze Applicate, Universit\`a di Roma ``La Sapienza'', 2004, to appear in Constr. Approx.
\bibitem{FJ} {M. Frazier, B. Jawerth}, {Decomposition of Besov spaces}, {Indiana Univ. Math. J.}, {34}, 777--799, 1985.
\bibitem{FJ2} {M. Frazier, B. Jawerth}, {Littlewood-Paley theory and the study of function spaces}, AMS, Providence RI, 1991.
\bibitem{Fuehr} H. F\"uhr, {\it The Abstract Harmonic Analysis of Continuous Wavelet Transforms}, Thesis, Technical University Munich, 2002.
\bibitem{GH} J-P. Gabardo,  D. Han, Frames associated with measurable spaces, Adv. Comput. Math. 18, No.2-4, 127--147, 2003.
\bibitem{G} P. Gr\"obner, {\it Banachr\"aume Glatter Funktionen und Zerlegungs-Methoden}, Ph.D. thesis,
University of Vienna
, 1992.
\bibitem{gro} K. Gr{\"o}chenig, Describing functions: atomic decompositions versus frames, Monatsh. Math., 112, 1--41, 1991.
\bibitem{Gr} {K. Gr\"ochenig}, {\it Foundations of Time-Frequency Analysis}, Birkh\"auser Verlag, 2001.
\bibitem{gro1} {K. Gr\"ochenig}, Localization of frames, Banach frames, and the invertibility of the frame operator,  J. Four. Anal. Appl., 10(2), 2004. 
\bibitem{gl} K. Gr\"ochenig, M. Leinert, Symmetry of matrix algebras and symbolic calculus for infinite matrices, Preprint, 2003.
\bibitem{GM} A. Grossmann, J. Morlet, and T. Paul, Transforms associated to square integrable group representation, I, J. Math. Phys., 26(10), 2473--2479, 1985.
\bibitem{GM1} A. Grossmann, J. Morlet, and T. Paul, Transforms associated to square integrable group representation, II, Examples, Ann. Inst. H. Poincar\'e, 25, 293--309, 1986.
\bibitem{j} S. Jaffard, Propri\'et\'es des matrices ``bien localis\'ees'' pre\`e de leur diagonale et quelques applications. Ann. Inst. H. Poicar\'e Anal. Non Lin\'eaire, 7(5), 461--476, 1990.
\bibitem{HL2} {J.A. Hogan, J.D. Lakey}, {Embeddings and uncertainty principles for generalized modulation spaces}, 
Chapter 4 in ``Modern Sampling Theory: Mathematics and Applications'', J. J. Benedetto and P. J. S. G. Ferreira (Eds.), Birkh\"auser, Boston, 73--105, 2000.
\bibitem{AA} A. L. Hohou\'eto, T. Kengatharam, S. T. Ali, and J. P. Antoine, Coherent states lattices and square integrability of representations, Preprint, 2003.
\bibitem{HN} M. Holschneider, B. Nazaret, An interpolation family between Gabor and wavelet transformations. Application to differential calculus and construction of anisotropic Banach spaces, to appear in Adv. In Partial Diff. Eq., ``Nonlinear Hyperbolic Equations, Spectral Theory, and Wavelet Transformations'' (Albeverio, Demuth, Schrohe, Schulze Eds.), Wiley 2003.
\bibitem{K} G. Kaiser, {\it A Friendly Guide to Wavelets}, Boston: Birkh\"auser. xiv, 1994.
\bibitem{PS} L. P\"aiv\"arinta, E. Somersalo, A generalization of the Calderon-Vaillancourt theorem to $L^p$ and $h^p$, Math. Nachr. {138}, 145--156, 1988.
\bibitem{P} J. Peetre, Paracommutators and minimal spaces, in ``Operators and Function Theory'' (S. C. Powers, Ed.), NATO ASI Series, Reidel, Dordrecht, 1985.
\bibitem{hr} H. Rauhut, {\it Time-Frequency and Wavelet Analysis of Functions with Symmetry Properties}, Ph.D. thesis, 
Technical University Munich, 2004.
\bibitem{hr1} H. Rauhut, Wavelet transforms associated to group representations and functions invariant under symmetry
groups, Preprint 2003, to appear in Int. J. Wavelets, Multisc. and Inf. Proc.
\bibitem{hr2} H. Rauhut, Banach frames in coorbit spaces consisting of elements which are invariant under symmetry groups, Preprint 2003.
\bibitem{Roes} M. R\"osler, {Radial wavelets and Bessel-Kingman hypergroups}, Preprint, 1997.
\bibitem{T1} {B. Torresani}, {Wavelets associated with representations of the Affine Weyl-Heisenberg group}, J. Math. Phys., {32}, 1273--1279, 1991.
\bibitem{T2} {B. Torresani}, {Time-frequency representation, wavelet packets and optimal decomposition}, Ann. Inst. H. Poincar\'e, {56}, 215--234, 1992.
\bibitem{T3} {B. Torresani}, Position-frequency analysis for signals defined on spheres, Signal Process., 43(3), 341--346, 1995.
\bibitem{Tr1} {H. Triebel}, {\it Theory of Function Spaces}, 
Birkh\"auser, Basel, 1983.
\bibitem{Tr2} {H. Triebel}, {\it Theory of Function Spaces II}, Birkh\"auser, 1992.
\end{thebibliography}
\end{document}